\documentclass{siamltex}
\usepackage{amssymb}
\usepackage{amsmath}
\usepackage{amsfonts}
\usepackage{graphicx}
\usepackage{xcolor}
\usepackage{float}
\usepackage{multicol}
\usepackage{subcaption}
\usepackage{multirow}
\newtheorem {remark}[theorem]{Remark}

\title{Modular data assimilation for flow prediction}
\author{Aytekin \c{C}{\i}b{\i}k \thanks{Department of Mathematics, Gazi University, Ankara, 06550, Türkiye (abayram@gazi.edu.tr)} \and Rui Fang\thanks{Department of Mathematics, University of Pittsburgh, Pittsburgh, PA, 15260, USA (ruf10@pitt.edu)} \and William Layton\thanks{Department of Mathematics, University of Pittsburgh, Pittsburgh, PA, 15260, USA (wjl@pitt.edu). Corresponding author: William Layton. The work of William Layton and Rui Fang was partially supported by NSF grant DMS 2410893. The work of Aytekin \c{C}{\i}b{\i}k was partially supported by TUBITAK with the BIDEB-2219 grant.}}
\date{11 October 1999}

\begin{document}
\maketitle
\begin{abstract}
This report develops several modular, 2-step realizations (inspired by Kalman filter algorithms) of nudging-based data assimilation 
\begin{equation*}
\begin{array}{cc}
Step{{{\text { }}}}1 & 
\begin{array}{c}
\frac{\widetilde {v}^{n+1}-v^{n}}{k}+v^{n}\cdot \nabla \widetilde {v}%
^{n+1}-\nu \triangle \widetilde {v}^{n+1}+\nabla q^{n+1}=f(x){{{{\text { }}}}%
} \\ 
\nabla \cdot \widetilde {v}^{n+1}=0%
\end{array}
\\ 
Step{{{\text { }}}}2 & \frac{v^{n+1}-\widetilde {v}^{n+1}}{k}-\chi
I_{H}(u(t^{n+1})-v^{n+1})=0.%
\end{array}%
\end{equation*}
Several variants of this algorithm are developed. Three main results are developed. The first is that if $I_{H}^{2}=I_{H}$, then Step 2 can be rewritten as the explicit step 
\begin{equation*}
v^{n+1}=\widetilde {v}^{n+1}+\frac{k\chi }{1+k\chi }[I_{H}u(t^{n+1})-I_{H}%
\widetilde {v}^{n+1}].
\end{equation*}
This means Step 2 has the greater stability of an implicit update and the lesser complexity of an explicit analysis step. The second is that the basic result of nudging (that for $H$ \textit{small enough} and $\chi $\ \textit{large enough} predictability horizons are infinite) holds for one variant of the modular algorithm. The third is that, for \textit{any} $H>0$ and \textit{any} $\chi>0$, one step of the modular algorithm decreases the next step's error and \textit{increases} (an estimate of) predictability horizons. A method synthesizing assimilation with eddy viscosity models of turbulence is also presented. Numerical tests are given, confirming the effectiveness of the modular assimilation algorithm. The conclusion is that the modular, 2-step method overcomes many algorithmic inadequacies of standard nudging methods and retains a robust mathematical foundation.
\end{abstract}


\section{Introduction}
Predicting the future state of a flow, here the internal $2d$ or $3d$ flow of an incompressible viscous fluid in a domain $\Omega$, 
\begin{align}
u_{t}+u\cdot \nabla u-\nu \triangle u+\nabla p& =f(x){{{{\text{, and }}}}}%
\nabla \cdot u=0,{{{{\text{ in }}}}}\Omega , 0<t\leq T, \\
u& =0{{{{\text{ on }}}}}\partial \Omega {{{{\text{ and }}}}}u(x,0)=u_{0}(x),
\end{align}
is a continuing challenge with impact in science and technology. Perfect knowledge of flow data, especially initial conditions, is not possible. Exponential separation of velocity trajectories causes initial errors to grow exponentially until errors saturate at ${{{\mathcal{O}}}}(1)$ levels, inspiring a finite (and small) predictability horizon. Lorenz \cite{L85} noted that at this point the entire flow simulation has become, in effect, an expensive random number generator. The growth of these errors is slowed (and the predictability horizon extended), Kalnay \cite{K03}, by incorporating / assimilating observations of the flow 
\begin{equation*}
u_{obs}(x,t)=I_{H}u(x,t),{\text{ }}H\simeq {\text{ observation spacing}}.
\end{equation*}
Current best practice employs Kalman filter variants, an explicit, 2-step process, for data assimilation. \textit{Practice} thus deviates from \textit{theory}, where nudging / Newtonian damping / time relaxation (see Breckling, Neda, Hill \cite{BNH17} for connections) has developed a nearly complete theoretical foundation since the 2014 work of Azouani, Olson and Titi \cite{AOT14}. However, nudging still faces algorithmic issues, Section 1.1, that maybe one reason for a preference for a Kalman filter variant. Herein, inspired by the computational efficiency of Kalman filter algorithms, we develop and test modular nudging algorithms that ameliorate or even eliminate practical difficulties with nudging-based
assimilation. Since these are different methods than standard nudging, they require an analytic foundation, also developed herein, to understand their reliability.

Suppressing temporarily space discretization (which has a secondary role) to present the ideas, the first modular, 2-step method follows. With boundary and initial conditions $v=0{{{{\text{ on }}}}}\partial \Omega ,$ ${{{{\text{and }}}}}v(x,0)=v_{0}(x)$, with initial errors $v_0(x)\neq u_0(x)$, pick $\chi >0$ and proceed by forecast and analysis steps 
\begin{equation}
\begin{array}{cc}
Step{{{\text{ }}}}1 & 
\begin{array}{c}
\frac{\widetilde{v}^{n+1}-v^{n}}{k}+v^{n}\cdot \nabla \widetilde{v}%
^{n+1}-\nu \triangle \widetilde{v}^{n+1}+\nabla q^{n+1}=f(x){{{{\text{ }}}}}
\\ 
\nabla \cdot \widetilde{v}^{n+1}=0,{{{{\text{ }}}}}%
\end{array}
\\ 
Step{{{\text{ }}}}2A & \frac{v^{n+1}-\widetilde{v}^{n+1}}{k}-\chi
I_{H}(u(t^{n+1})-v^{n+1})=0.
\end{array}
\label{eq:Nudging2Step}
\end{equation}
This combination is interpretable a type of operator split, 1-step nudged approximation, (\ref{2step_B}) below. The form of (\ref{eq:Nudging2Step}) resolves algorithmic limitations, detailed in Section 1.1 and Section 7 below, of standard nudging. While open problems remain (Section 7), this combination promises an assimilation method with less computational complexity than Kalman filters and a mathematical foundation approaching standard nudging. Various choices of $I_{H}$ have been studied. Herein, we make the common choice of taking $I_{H}$ to be interpolation, restriction, or projection (associated with a coarse mesh scale $H$) into a finite-dimensional space $X^{H}$. Many of these choices satisfying $I_{H}=I_{H}^{2}$.

In Sections 4 and 5, to shorten the technical error analysis, $I_{H}$ will there be a projection with respect to either the $L^{2}$ inner product (denoted as usual $(\cdot, \cdot)$) or $H^{1}$ semi-inner product. Section 1.4 briefly considers modular spectral nudging where $I_{H}$ is a local averaging: 
\begin{eqnarray*}
(u-I_{H}u,w^{H}) &=&0{{{\text{ or }}}}(\nabla \left[ u-I_{H}u\right] ,\nabla
w^{H})=0,\forall {{{\text{ }}}}w^{H}\in X^{H}\text{ (Sections 4,5),} \\
I_{H} &=&\text{ a local averaging operator (1.4).}
\end{eqnarray*}
Step 1 can be replaced by one step of any standard method for the (laminar) Navier-Stokes equations (NSE). Steps 1 and 2 together resemble (various) Kalman filter realizations where $\chi$ would be replaced by an approximate error covariance matrix and Step 2 would be made explicit by replacing $I_{H}(u(t^{n+1})-v^{n+1})$ by $I_{H}(u(t^{n+1})-\widetilde{v}^{n+1})$. The explicit treatment of $I_Hv$ in Step 2 is one reason for the greater computational efficiency of Kalman filter variants.

For (\ref{eq:Nudging2Step}) and variants, we present 3 main results. The first, Section 2, is that \textit{Step 2A does not require assembly and solution of a linear system
coupling fine and coarse scales.} If $I_H^2 =I_H$, it can be executed \textit{explicitly} by 
\begin{equation}
v^{n+1}=\widetilde{v}^{n+1}+\frac{k\chi }{1+k\chi }[I_{H}u(t^{n+1})-I_{H}
\widetilde {v}^{n+1}].  \label{eq:ExplicitUpdateIntro}
\end{equation}
This formula gives the same complexity, without the stability problems of an explicit analysis step, as when Step 2 is made explicit by replacing $I_{H}(u(t^{n+1})-v^{n+1})$ by $I_{H}(u(t^{n+1})-\widetilde{v}^{n+1})$. This explicit update is possibly only for modular, 2-step assimilation and not for standard nudging. The second main result mirrors the striking property of the usual nudged approximation that for any Reynolds number, for $H$ \textit{small enough} and $\chi$\ \textit{large enough} the resulting approximation is accurate \textit{uniformly in time.}
The \textquotedblleft \textit{how small and large enough}\textquotedblright\ restrictions required are comparable to ones for the classical nudged approximation in, e.g.,\cite{AOT14}, \cite{BP21}, \cite{CGJP22}, \cite{LRZ19}.
The proof, Section 4, requires new ideas and is more intricate than for standard nudging of Section 1.1. An estimate of predictability horizons can be obtained from data at $t_n$ and $t_{n+1}$. The third main result, Section 3, addressed the case when the (unrealistic \cite{CFLS25}) \textquotedblleft $H$ \textit{%
small enough} and $\chi$\ \textit{large enough}\textquotedblright conditions fail. Section 3 shows that for \textit{any} $H>0,\chi >0$, the method (\ref{eq:Nudging2Step}) decreases error and increases this estimate of the predictability horizon (defined in Section 3). Since data for infinite time is not available and the density of data locations can be inflexible, this (simply proven) result is significant.

Naturally, the most interesting case is where assimilation is applied for higher Reynolds number flows. In this case, a turbulence model would also normally be included in the system. The form of the 2-step method and previous work in \cite{LRT12} suggest that Step 2 can be augmented to implement assimilation and an eddy viscosity turbulence model. Several realizations of this combination are possible. We investigate in Section 5
the alternate Step 2: 
\begin{equation}\label{2step_B}
\begin{array}{cc}
Step{{{\text{ }}}}1 & 
\begin{array}{c}
\frac{\widetilde{v}^{n+1}-v^{n}}{k}+v^{n}\cdot \nabla \widetilde{v}%
^{n+1}-\nu \triangle \widetilde{v}^{n+1}+\nabla q^{n+1}=f(x){{{{\text{ }}}}}
\\ 
\nabla \cdot \widetilde{v}^{n+1}=0,{{{{\text{ }}}}}%
\end{array}
\\ 
Step{{{\text{ }}}}2B & 
\begin{array}{c}
\frac{v^{n+1}-\widetilde{v}^{n+1}}{k}-\chi I_{H}(u(t^{n+1})-v^{n+1}) \\ 
-\nu \triangle \left( v^{n+1}-\widetilde{v}^{n+1}\right) -\nabla \cdot (\nu
_{turb}(\widetilde{v}^{n+1})\nabla v^{n+1})=0.%
\end{array}%
\end{array}%
\end{equation}
We show in Section 5 that the additional term $-\nu \triangle \left( v^{n+1}-\widetilde{v}^{n+1}\right) $\ (with, as a first step, $\nu _{turb}(\widetilde{v}^{n+1})
=0$) does not alter the essential property of the method having an infinite predictability horizon under conditions on $\chi, H$. Testing and analysis with $\nu _{turb}(\widetilde{v}^{n+1})\neq 0$ is an open problem.

\begin{remark}
Steps 1 and 2 are close to an operator splitting of (\ref{eq:StandardNudged}) below so concern of this combination possibly inducing an accuracy barrier is reasonable. There are 3 cases where operator splitting does not induce an accuracy barrier: when the term split is small (or higher order), for special solutions detailed in Lanser and Verwer \cite{LV99} and when the term is exactly zero for the true solution. This last is the case here since $\chi I_{H}(u(t^{n+1})-v^{n+1})=0$ when $v=u$. 
\end{remark}

\subsection{Standard Nudging}

We again (temporarily) suppress space discretization and use a simple time discretization. The approximation $v^{n+1}(x)$ to $u(x,t^{n+1})$\ is: ${{{{%
\text {select parameter }}}}}\chi ,{{{{\text { set\ }}}}}v^{0}=v_{0}(x),v=0$
on $\partial \Omega $,${{{{\text { and solve}}}}}$ 
\begin{gather}
\frac{v^{n+1}-v^{n}}{k}+v^{n}\cdot \nabla v^{n+1}-\nu \triangle v^{n+1}
\label{eq:StandardNudged} \\
-\chi I_{H}(u(t^{n+1})-v^{n+1})+\nabla q^{n+1}=f(x){{{{\text {, }}}}}\nabla
\cdot v^{n+1}=0.  \notag
\end{gather}
There are several issues with the practical use of standard (1 step) nudging calling for algorithmic improvements.

\textbf{Legacy codes:} Nudging in (1.6) requires an intrusive modification of existing computational fluid dynamics (CFD) codes.

\textbf{Coefficient matrix fill-in:} If (\ref{eq:StandardNudged}) is discretized on a spatial mesh of size $h \ll H$, the term $\chi I_{H}v^{n+1}$ couples variables across the observational macro-elements, increasing fill-in and operations (memory and complexity) dramatically. There are cases when this can be ameliorated in Rebholz and Zerfas \cite{Rebholz_Zerfas_2021}. These require specific choices excluding some realizations. 

\textbf{Mesh communication:} Nudging also requires communication between the 
\newline
$H-$mesh and the $h-$mesh. Its availability depends upon the software platform used.

\textbf{Conditioning:} When $\chi \gg 1$, the nudging term increases the condition number of the associated coefficient matrix by adding a term with $\chi$ as a multiplier. Similar issues occur here when (1.4) is not used, Section 2.2. We note that the system from Step 2A will be symmetric positive definite (SPD) and thus less sensitive to high condition numbers than (1.6). 

\textbf{The }$H, \chi $\textbf{\ conditions:} The precise form of $``H$ \textit{small enough} and $ \chi$ \textit{large enough" } in, e.g. 3d is${{{{%
\text{ }}}}}\chi \gtrsim {{{{\mathcal{O}}}}}({{{{\mathcal{R}}}}}e^{5}){{{{%
\text{ and \ }}}}}H\lesssim {{{{\mathcal{O}}}}}({{{{\mathcal{R}}}}}e^{-3})$,
smaller than the Kolmogorov micro-scale, \cite{CFLS25}. Also, with $\chi
\gtrsim {{{{\mathcal{O}}}}}({{{{\mathcal{R}}}}}e^{5})$\ the nudged model no
longer \textquotedblleft \textit{transmits information from data-rich to data-poor regions}\textquotedblright, p.16 \cite{K03}, and \textquotedblleft \textit{prevents buildup of small scale variability} \textquotedblright , \cite{ODD12b}. Estimates of good $\chi$ values obtained from computational test vary but $\chi \sim 40\%$ of the predictability horizon $\left( \ll{{{{\mathcal{R}}}}}e^{5}\right) $ in \cite{ODD12a} is a typical example.

\textbf{The user-chosen parameter }$\chi$\textbf{ must be chosen} with some progress but still not definitive analytic guidance.

\subsection{Related work}

The goal of nudging based data assimilation, from 1964 Luenberger \cite{L64}, is to combine (noisy) data with (approximate) deterministic laws to produce better forecasts than either alone, Kalnay \cite{K03}. The implementation of standard nudging in \textit{research} codes is straightforward if communication between meshes is not required, such as when $v^{n+1}, \widetilde {v}^{n+1}$ are piecewise quadratic and $I_{H}v^{n+1}$ is piecewise constant \textit{on the same mesh}. On the other hand, in more general (and useful) cases, communication between meshes can be challenging in research codes and impossible in production codes. The theoretical analysis of standard nudging is quite general, highly developed and shows the great potential of data assimilation, e.g., Azouani, Olson and Titi \cite{AOT14}, Biswas and Price \cite{BP21},
Cao, Giorgini, Jolly and Pakzad \cite{CGJP22}, Larios, Rebholz, and Zerfas \cite{LRZ19} (among hundreds of papers). In particular, the now classical proof that for any Reynolds number and suitable $\chi, H$ the nudged approximation can have an infinite predictability horizon is a spectacular theoretical accomplishment. For applications see, for example, Lakshmivarahan and Lewis \cite{LL13}, Kalnay \cite{K03}, Navon \cite{N98}, and for interesting and important recent work in various directions Carlson, Hudson, Larios, Martinez, Ng, Whitehead \cite{CHLMN21} and Carlson, Farhat, Martinez, Victor \cite{CFMV24}.

\subsection{Preliminaries}

The notation for the $L^{2}(\Omega )$ norm, and inner product is $\Vert \cdot \Vert $ and $(\cdot , \cdot)$, respectively. Using $\Vert \cdot \Vert_{L^{p}}$, we indicate the $L^{p}(\Omega )$ norm. The numerical tests use a standard finite element discretization, see \cite{L08} for details. We also use the following inequalities. The first estimate of this type was derived in 1958 by Nash \cite{N58}. These are proven in Lemmas 1, 2, 3 pages 10, 11, 12, Section 1.1 of Ladyzhenskaya \cite{L69} with constant, $2^{1/4}$, $2^{1/2}$, $48^{1/6}$ respectively. Soon thereafter work of Gagliardo \cite{gagliardo1959ulteriori} and Nirenberg \cite{nirenberg1959elliptic} independently derived and extended the bounds. Others improved the second and the third constants to those below in particular the $2 \cdot 3^{-1/2}$ is improvable to $(3 \pi)^{-1/2}$. 

\begin{lemma}\label{norms_property} For any vector function $u:{{{{\mathbb{R}}}}}%
^{d}\rightarrow {{{{\mathbb{R}}}}}^{d}$ with compact support and with the
indicated $L^{p}$ norms finite, 
\begin{align*}
\Vert u\Vert _{L^{4}({{{{\mathbb{R}}}}}^{2})}& \leq 2^{1/4}\Vert u\Vert^{1/2}\Vert \nabla u\Vert^{1/2},\ (d=2), \\
\newline
\Vert u\Vert _{L^{4}({{{{\mathbb{R}}}}}^{3})}& \leq \frac{4}{3\sqrt{3}}\Vert
u\Vert ^{1/4}\Vert \nabla u\Vert ^{3/4},\ (d=3),\newline
\\
\Vert u\Vert _{L^{6}({{{{\mathbb{R}}}}}^{3})}& \leq \frac{2}{\sqrt{3}}\Vert
\nabla u\Vert ,\ (d=3).
\end{align*}%
\end{lemma}
For the nonlinear term, we use the following bounds \cite{L08}, 
\begin{equation}
\begin{split}
(u\cdot \nabla v,w)& \leq C_{2}\Vert \nabla u\Vert \Vert \nabla v\Vert \Vert
\nabla w\Vert , \\
& \leq C_{2}\sqrt{\Vert u\Vert }\sqrt{\Vert \nabla u\Vert }\Vert \nabla
v\Vert \Vert \nabla w\Vert,
\end{split}%
\end{equation}%
for $u,v,w\in (H_{0}^{1})^d$, $d=2$ or $3$.

\subsection{Spectral nudging}

In spectral nudging, the nudging term is tuned to only act on the solution's large scales, e.g., Omrani, Drobinski, Dubos \cite{ODD12a} \cite{ODD12b}, Radu, D\'{e}qu\'{e}, Somot \cite{RDS08}, like time relaxation algorithms in large eddy simulation, analyzed in \cite{LN07}. The essential question is how to do this effectively when all the data is discrete. As a first analysis step, we study here (non-discrete) \textit{continuum} averaging by a Germano differential filter \cite{Germano86} as a reasonable and simplifying proxy for local averaging of discrete data.

\begin{definition}
Given a velocity $w \in H_{0}^{1}(\Omega)^{d},$ its average over scale $H$, $\overline{w}=I_{H}w$ is the solution of
\begin{equation}\label{average_operator}
-H^{2}\triangle \overline{w}+\overline{w}=w\text{ in }\Omega \text{\ and }%
\overline{w}=w(=0)\text{\ on }\Omega.
\end{equation}
\end{definition}

The spectral nudging realization we study is then $\overline{w}=I_{H}w$ in Step 2
\begin{equation*}
Step\text{ }2A:\text{ \ }\frac{v^{n+1}-\widetilde{v}^{n+1}}{k}-\chi (%
\overline{u(t^{n+1})}-\overline{v^{n+1}})=0.
\end{equation*}
The analysis in Sections 4, 5 shows that nudging's basic result of an infinite predictability horizon for $H$ small and $\chi$ large will hold for the modular algorithm (1.3) and variants provided $\overline{w}=I_{H}w$ satisfies 3 properties. The result that nudging increases predictability horizons in Section 3 depends on a fourth
property. This averaging (\ref{average_operator}) is an elliptic-elliptic singular perturbation problem for which many estimates, global and local (e.g., \cite{L12}) are known. The first three essential properties are proven in
Appendix A of \cite{LR12}. The fourth follows from $\left( -H^{2}\triangle +1\right)^{-1}$ being an SPD operator. These four are summarized in the next proposition.

\begin{proposition}
For $w\in H_{0}^{1}(\Omega )^{d}$ we have 
\begin{eqnarray*}
\|I_{H}w\| &=&\|\overline{w}\|\leq \|w\|,\text{ \ }\|\nabla
I_{H}w\|=\|\nabla \overline{w}\|\leq \|\nabla w\|, \\
\|w-\overline{w}\| &\leq &\frac{1}{2}H\|\nabla w\|\text{ and }(w,\overline{w}%
)>0\text{ for }w\neq 0.
\end{eqnarray*}
\end{proposition}

The error estimates established for nudging with $I_{H}$ a projection follow for this realization of spectral nudging by the same arguments. The result that Step 2 can be implemented as an explicit step does not follow however since $I_H^2 \neq I_H$ for local averaging.

\section{Step 2 can be explicit}

Executing Step 2 (in finite element methods with $I_{H}$ a projection) generally involves discretization in space (on a mesh of width $h\ll H$) the equation 
\begin{gather*}
\frac{v^{n+1}-\widetilde {v}^{n+1}}{k}-\chi I_{H}(u(t^{n+1})-v^{n+1})=0,{{%
\text { }}}or \\
v^{n+1}+k\chi I_{H}v^{n+1}=\widetilde {v}^{n+1}+k\chi I_{H}u(t^{n+1}){{%
\text { }}}
\end{gather*}
then solving the associated linear system for the discrete approximation to $v^{n+1}$. For this subsection, we suppress the superscript \textquotedblleft 
$n+1$\textquotedblright and focus on dependence on $h$ and $H$. We consider 2 approaches to this step.

The first, Section 2.1 is based on a remarkable (but easily proven) observation that \textit{Step 2 can be executed explicitly and requires only 1 explicit application of the interpolation / restriction projection operator $I_H$}. The second is by direct implementation of Step 2 in a finite element software package. This approach involves assembly and solution of a linear system for $v^{n+1}$'s degrees of freedom. Since $\chi \gg 1$ (and current theory considers $\chi \to \infty$, for example, \cite{Diegel2025}) we develop an estimate of the condition number of this linear system as a function of $\chi$ in Section 2.2.  

\subsection{\textbf{Approach 1: An explicit calculation} }
The main result here is that if $I_H^2=I_H$, $v^{n+1}$ can be calculated by an explicit update
\begin{equation}\label{form_a}
    v^{n+1} = \widetilde{v}^{n+1} +\frac{k\chi}{1+k\chi} I_H(u^{n+1}-\widetilde{v}^{n+1}).
\end{equation}
This implies that (if so implemented) modular nudging has lower complexity than an explicit Kalman filter step. 

\begin{proposition}\label{prop_equv_form} In all cases Step 2A satisfies
\begin{equation}\label{form_b}
    v^{n+1} = \widetilde{v}^{n+1} +\frac{k\chi}{1+k\chi} I_H(u^{n+1}-\widetilde{v}^{n+1}) + \frac{\left(k\chi\right)^2}{1+k\chi}\left[I_H -I_H^2\right] (u^{n+1}-v^{n+1}).
\end{equation}
Thus, if $I_H=I_H^2$, (\ref{form_a}) above holds.
\end{proposition}
\begin{proof}
The relations (\ref{form_a}), (\ref{form_b}) follow by, essentially, rearrangement of Step 2. Dropping the superscripts \textquotedblleft n+1\textquotedblright for clarity, we have
\begin{equation}
    v=-k\chi I_Hv + \widetilde{v} + k \chi I_H u.
\end{equation}
Applying $I_H$ to the above gives 
\begin{equation}
    I_H v = \frac{1}{1+k\chi} I_H\widetilde{v}+\frac{k \chi}{1+k\chi} I_H u + \frac{k\chi}{1+k\chi} \left[ I_H^2 -I_H\right] u -\frac{k\chi}{1+k\chi} \left[I_H^2-I_H \right] v.
\end{equation}
Substituting the right-hand side (RHS) of the above for $I_H v$ in the term \textquotedblleft $-k\chi I_H v$\textquotedblright of the first equation and rearranging completes the proof.
\end{proof}

\begin{remark} On $I_H^2=I_H$. The condition $I_H^2=I_H$ holds trivially if $I_H$ is an interpolation of point values at locations that are nodes of a mesh. If $I_H$ is a local averaging operator, Section 1.4, then it fails as double averaging expands the effective averaging radius.  If Step 2 is considered without space discretization, $I_H^2=I_H$ holds for any projection. With space discretization considered, the condition becomes more complex. 
\end{remark}

\begin{remark} On $I_H^2 = I_H$ with space discretization. Suppose Step 2 is implemented by a finite element space discretization. Then $v^{n+1}, \widetilde{v}^{n+1}\in X^h \subset (H^1_0(\Omega))^d$ with $X^h$ a finite element space based on a mesh width $h$ ($h\ll H$). Let
\begin{equation}
    I_h^0 = L^2 \text{ projection into } X^h.
\end{equation}
Then the finite element space discretization of Step 2 can be written as follows: find $v^{n+1} \in X^h$ satisfying 
\begin{equation}
    \frac{v^{n+1}-\widetilde{v}^{n+1}}{k} = \chi I_h^0 I_H(u^{n+1}-v^{n+1}).
\end{equation}
Thus, Proposition \ref{prop_equv_form} applies after space discretization if
$(I_h^0 I_H)^2=I_h^0 I_H$, i.e. if $I_h^0 I_H$ is a projection. 

Suppose $I_H$ is an $L^2$ Projection into some coarse mesh subspace $X^H$. A classic result in linear algebra is that $I_h^0 I_H$ is a projection if and only if $I_h^0$ and $I_H$ commute. It follows then that $I_h^0 I_H$ is a projection provided $X^H \subset X^h$. On the other hand, a common choice for academic tests (because it does not require communication between 2 meshes) is, on the same mesh,
\begin{equation}
\begin{gathered}
    X^h = C^0 \text{ piecewise quadratics},\\
    X^H = \text{ discontinuous piecewise constants}.
    \end{gathered}
\end{equation}
With this choice $I_h^0 I_H \neq I_H I_h^0$ and (\ref{form_a}) fails. The analysis of this case must account for the third term as the RHS of (\ref{form_b}), either in the algorithmic realization or as a consistency error term. For the latter we note a classic result of Crouzeix and Thom{\'e}e \cite{CT_87} that for $L^2$ projections into finite element spaces (under mild conditions)
\begin{equation}
    \|(I_H - I_H^2)e\| = \|(I-I_H)(I_He)\|\leq C_1 H \|\nabla I_H e\|\leq C_2 H \|\nabla e\|.
\end{equation}
Finally, we note that generic non-commutativity of projections is the basis upon which convergence of alternating projection methods is built, beginning with work of Kaczmarz \cite{SK1937} and von Neumann \cite{VN49} and continuing actively, e.g. Escalante and Raydan \cite{ER11}. 
\end{remark}

\subsection{\textbf{Approach 2:} \textbf{Assemble and solve the linear system for }$%
v^{n+1}$} With some finite element method (FEM) software, the quickest (in programmer effort) way to get results is just to assemble and solve, typically by conjugate gradient (CG) methods, the system of Step 2. This can be efficient in computational complexity became the resulting coefficient matrix is SPD if $I_{H}$ is an $L^{2}$ projection. This approach has two concerns: \textit{conditioning} (since $\chi \gg 1$) and \textit{complexity} (since $H\gg h$).

\textbf{Conditioning in Approach 2.} The analysis of condition numbers next follows standard steps. To shorten the presentation, we assume $I_{H}$ is an $L^{2}$
projection. Denote the velocity space $X:=\left(H_{0}^{1}(\Omega)\right)^{d},$ $d=2, 3$. Finite element velocity sub-spaces associated with fine and coarse meshes are denoted $X^{h}\subset X$ and $X^{H}\subset
L^{2}(\Omega )^{d}$, see e.g. \cite{L08} for details. After discretization in space, Step 2 is: Given $\widetilde {v}^{h}\in X^{h}$, find $v^{h}\in X^{h}$, satisfying, for $\phi ^{h}\in X^{h}$
\begin{equation*}
(v^{h},\phi ^{h})+k\chi (I_{H}v^{h},\phi ^{h})=(\widetilde {v}^{h},\phi^{h})+k\chi (I_{H}u,\phi ^{h}).
\end{equation*}
Let $\phi _{i}^{h},i=1,\cdot \cdot \cdot ,N_{h}$ be a basis for $X^{h}$ and $\phi _{i}^{H}$ a similar basis for $X^{H}$. Expanding 
\begin{equation*}
v^{h}=\sum _{i=1}^{N_{h}}c_{i}\phi _{i}^{h},
\end{equation*}
the system for the vector $c$ of coefficients is then 
\begin{equation*}
\sum _{i=1}^{N_{h}}c_{i}\left (\left [1+k\chi I_{H}\right ]\phi
_{i}^{h},\phi _{j}^{h}\right )=(f,\phi _{j}^{h}).
\end{equation*}
Let $M^{s}$ denote the finite element mass matrix $M_{ij}^{s}=($ $\phi_{i}^{s},$ $\phi _{i}^{s}),s=h$ $\&$ $H$. The coefficient matrix is thus, $%
i,j=1,\cdot \cdot \cdot ,N_{h},$ 
\begin{equation*}
A_{ij}=\left (\left [1+k\chi I_{H}\right ]\phi _{i}^{h},\phi _{j}^{h}\right
)=\left (\phi _{i}^{h},\phi _{j}^{h}\right )+k\chi \left (I_{H}\phi
_{i}^{h},I_{H}\phi _{j}^{h}\right ).
\end{equation*}
We first prove that, under mild conditions, $cond(A)\leq C(1+k\chi)$, as expected.

\begin{theorem}
Assume $I_{H}$ is an $L^{2}$ projection. Assume that for $s=h$ $\&$ $H$\ and all $v\in X^{s},$ $v=\sum_{i=1}^{N_{s}}c_{i}\phi _{i}^{s}(x),$ $\|v\|$ and $%
\sqrt{N_{s}^{-1}\sum_{i=1}^{N}c_{i}^{2}}$ are uniform-in-$s$ equivalent norms. Then, 
$$cond(A)\leq C(1+k\chi).$$   
\end{theorem}
\begin{proof}
For any $N_{s}$\ vector $c$, define $w^{s}=\sum_{i=1}^{N_{s}}c_{i}\phi
_{i}^{s}(x) \in X^{s}.$ Norm equivalence implies $\|w\|$ and $\sqrt{N_{s}^{-1}\sum _{i=1}^{N_{s}}c_{i}^{2}}$ are uniformly in $s=h$ $\&$ $H$ equivalent norms. Since $A$ is SPD, we estimate $\lambda (A)$. Let $\lambda, c$\ be an eigenvalue, eigenvector pair of $A$. Then, with $w=\sum
_{i=1}^{N_{h}}c_{i}\phi _{i}^{h}(x)$ 
\begin{equation*}
\lambda =\frac{c^{T}Ac}{c^{T}c}=\frac{\|w\|^{2}+k\chi \|I_{H}w\|^{2}}{|c|^{2}
}.
\end{equation*}
By norm equivalence of $\|w\|^{2}$ with $N_{h}^{-1}|c|^{2}$ and $0\leq \|I_{H}w\|^{2}\leq \|w\|^{2}$ we thus have $0 < CN_{h}\leq \lambda \leq CN_{h}(1+k\chi )$ and thus, uniformly-in $h\&H$, 
\begin{equation*}
cond(A)\leq C(1+k\chi).
\end{equation*}
\end{proof}

\textbf{A remark on complexity of approach 2.} We continue to assume $I_{H}$ is an $L^{2}$ projection. Since the step 2 coefficient matrix is SPD with condition number ${{\mathit{O}}}(1+k\chi)$, the natural solution procedure is by CG methods. At each step a residual must be calculated. There are two natural ways to compute
the residual: using linear algebra directly and by a combination of linear algebra and FEM assembly.

We first consider the direct approach via linear algebraic operations, easily implemented in FEM software packages that allow communication between 2 meshes. Introduce the variable $w^{H}:=I_{H}v^{h}$. Rewrite the system 
\begin{equation}
(v^{h},\phi ^{h})+k\chi (I_{H}v^{h},\phi ^{h})=(b^{h},\phi ^{h}){{\text {
where }}}(b^{h},\phi ^{h}):=(\widetilde {v}^{h},\phi ^{h})+k\chi
(I_{H}u,\phi ^{h})  \label{eqStepVarForm}
\end{equation}
as (for all $\phi ^{s},s=h$ $\&$ $H$ indicated) 
\begin{equation*}
(v^{h},\phi ^{h})+k\chi (w^{H},\phi ^{h})=(b^{h},\phi ^{h}){{\text { \ and \ 
}}}(w^{H},\phi ^{H})-(v^{h},\phi ^{H})=0.
\end{equation*}
Expand 
\begin{equation}
v^{h}=\sum _{i=1}^{N_{h}}c_{i}\phi _{i}^{h}(x),w^{H}=\sum
_{i=1}^{N_{H}}d_{i}\phi _{i}^{H}(x).  \label{eq:vhwHdefinition}
\end{equation}
Recall $M_{ij}^{s}=($ $\phi _{i}^{s},$ $\phi _{i}^{s}),s=h$ $\&$ $H$ and define the $N_{h}\times N_{H}$ matrix 
\begin{equation*}
B_{ij}:=(\phi _{i}^{h},\phi _{i}^{H}),1=1,\cdot \cdot ,N_{h},j=1,\cdot \cdot
\cdot ,N_{H}.
\end{equation*}
Substitute the two expansions into the above system. This leads to the coupled system for the coefficient vectors $c, d$ 
\begin{equation*}
\begin{bmatrix}
M^{h} & \mid & +k\chi B \\ 
- & + & - \\ 
-B^{T} & \mid & M^H
\end{bmatrix}
\begin{bmatrix}
c \\ 
- \\ 
d%
\end{bmatrix}
= 
\begin{bmatrix}
b \\ 
- \\ 
0%
\end{bmatrix}
.
\end{equation*}
Eliminating the vector $d$ gives the SPD reduced system coming from the formulation (%
\ref{eqStepVarForm}) 
\begin{equation*}
\left [M^{h}+k\chi B\left (M^{H}\right )^{-1}B^{T}\right ]c=b.
\end{equation*}
Its solution by CG methods is natural. Each step one must compute a residual 
\begin{equation*}
r=b-\left [M^{h}+k\chi B\left (M^{H}\right )^{-1}B^{T}\right ]c.
\end{equation*}
Introduce $d=$ $\left (M^{H}\right )^{-1}\left (B^{T}c\right ).$ The residual calculation becomes 
\begin{equation*}
r=b-M^{h}c-k\chi Bd{{\text { \ and \ }}}M^{H}d=B^{T}c.
\end{equation*}
The step $\left (M^{H}\right )^{-1}\left (B^{T}c\right )$\ in the residual calculation can be easily (since $cond(M^{H})={\mathcal{O}}(1))$ done by an inner iteration solving $M^{H}d=B^{T}c$. A second path to computing $\left(M^{H}\right)^{-1}\left (B^{T}c\right)$ in the residual calculation depends on the specified elements. If, for example, $X^{H}$ consists of piecewise constants on macro-elements (a common choice in tests), then $I_{H}v^{h}$ is the average of $v^{h}$ on macro-elements. This can be
computed by assembly using 1-point quadrature on the fine mesh elements comprising the macro-element. This approach depends on modifying the assembly routines in the code and is thus more intrusive.

\section{Modular nudging increases predictability horizons}

High Reynolds number flow is the motivating example of a chaotic system with a finite predictability horizon. The (induced) predictability horizon is fundamental to understanding fluid behavior. For example, Randall \cite{R12} p. 198 Ch.7 defines turbulence (imprecisely but evocatively) as \textquotedblleft \textit{A flow is turbulent if its predictability time is shorter than the time scale of interest}.\textquotedblright \ We now prove that adding Step 2 (the modular nudging step) after Step 1 increases (an estimates of) predictability horizons, calculated step by step. This requires a few preliminary definitions. First, we recall the definition given by Lorenz of a chaotic system (making precise the idea that chaotic means exactly but not approximately predictable).

\begin{definition}
A chaotic system is one with an absorbing ball and a positive Lyapunov exponent.
\end{definition}

For the NSE, the existence of an absorbing ball follows from the uniform boundedness of kinetic energy, \cite{L08}. Having a positive maximal Lyapunov exponent is observed in tests but analytical proof is elusive due to its formal, infinite time definition. We work instead with FTLE's (Finite Time Lyapunov Exponents), an idea developed by Boffetta, Giuliani, Paladin and Vulpiani \cite{BGPV98} and Cencini and Vulpiani \cite{CV13}. Errors are measured against the sought velocity, here $u(x,t)$.
\begin{definition}
Let $v(x,t)$ be an approximation to $u(x,t)$ and $\|\cdot \|_{\ast }$ a norm on the solution space. The FTLE calculated from the solution at times $T_1<T_2$
under the chosen norm is 
\begin{equation*}
\lambda (T_1, T_2) = \frac{1}{T_2-T_1}\ln \left(\frac{\|u(x,T_2)-v(x,T_2)\|_{\ast}}{
\|u(x,T_1)-v(x,T_1)\|_{\ast }}\right).
\end{equation*}
\end{definition}
The experimental doubling time of initial errors is then 
\begin{equation*}
\tau _{double}(T_1,T_2)=\frac{\ln 2}{\lambda (T_1, T_2)}>0.
\end{equation*}

When $\lambda (T_1,T_2)<0$, the error is shrinking. In that case $\frac{\ln(2)}{|\lambda|}$ is the error half life.
The $\varepsilon -$predictability horizon is 
\begin{equation*}
\tau ^{\varepsilon } (T_1, T_2)=\frac{1}{\lambda (T_1,T_2)}\ln \left (\frac{\varepsilon}{%
\|u(x,T_1)-v(x,T_1)\|_{\ast }}\right ).
\end{equation*}
Observe that \textit{decreasing $\lambda $\ increases the predictability horizon $\tau^{\varepsilon}.$ Since the natural logarithm is monotone increasing, decreasing the error at time $T_2$, $\|u(x,T_2)-v(x,T_2)\|$, decreases $\lambda$ and thus increases the predictability horizon.}

Define the errors 
\begin{equation*}
e^{n}(x)=u(x,t^{n})-v^{n}(x){{{\text { and }}}}\widetilde {e}%
^{n}(x)=u(x,t^{n})-\widetilde {v}^{n}(x).
\end{equation*}
Consider Step 2A. By subtraction the induced errors satisfy 
\begin{equation*}
\frac{e^{n+1}-\widetilde {e}^{n}}{k}+\chi I_{H}e^{n+1}=0.
\end{equation*}
If $I_{H}$ is an L$^{2}$ projection, take the inner product with $e^{n+1}$ and use the polarization identity. This yields 
\begin{equation*}
\frac{1}{2}\Vert e^{n+1}\Vert ^{2}-\frac{1}{2}\Vert \widetilde {e}%
^{n+1}\Vert ^{2}+\frac{1}{2}\Vert e^{n+1}-\widetilde {e}^{n+1}\Vert
^{2}+k\chi \Vert I_{H}e^{n+1}\Vert ^{2}=0.
\end{equation*}

\noindent This implies that $\|e^{n+1}\|<\|\widetilde{e}^{n+1}\|$ unless $v^{n+1}=\widetilde{v}^{n+1}$. Thus if we choose $T_1=t^{n}$, $T_2=t^{n+1}$, Step 2, having decreased the error at $t^{n+1}$, has increased the calculated predictability horizon compared to the predictability horizon calculated with Step 2 omitted.

\section{Infinite predictability horizons for $H$ small and $\chi$ large}

One striking result for nudging-based data assimilation is that for  $``H$ \textit{small enough} and $\chi$ \textit{large enough"} the nudged approximation has an infinite predictability horizon. Since the proof of this
is fundamental (and not difficult) for (\ref{eq:StandardNudged}), it is natural to investigate when it can be extended to the modular, 2-step method (\ref{eq:Nudging2Step}). At this point, we can extend this result to some but not all of the variants herein. Open questions are summarized in Section 7. We first present the stability of the two-step nudging with Step 2A. Then we derive an error estimate for the discrete-time, continuous-in-space of two step nudging with Step 2A.

We assume $X^H$ is a finite element space with mesh size $H$ such that $X^H \subset X =(H^1_0(\Omega))^d$, $d=2$ or $d=3$.

\begin{equation}\label{I_H_assump1}
\|(I-I_{H})w\|\leq C_{1}H\|\nabla w\|{{{\text{ for all }}}}w\in\left(
H_{0}^{1}(\Omega)\right)  ^{d},d=2{{{\text{ or }}}}3.
\end{equation}
We will use the following lemmas to establish the stability and error estimates.

We denote $I_H^0$ as the $L^2$ projection on $X^H$ and $I_H^1$ as the $H^1$ semi-inner product projection, which satisfy the following
\[
(u-I_{H}^0u,w^{H})=0{{\text{ and }}}(\nabla\left[  u-I_{H}^1u\right]  ,\nabla
w^{H})=0, {{\text{for }}} {{\text{all }}}w^{H}\in X^{H}.
\]
\begin{lemma} \label{layton2002}(Lemma 3, \cite{L02}) For all $w \in X^H$,
\begin{equation}
\nabla I_H^1 (w) = I_H^0 (\nabla w).
\end{equation}
\end{lemma}
Lemma \ref{grad_v_bound} shows that $\|\nabla v^{n+1}\|$ is bounded by terms involves $\nabla u(t_{n+1})$ and $\nabla \widetilde{v}^{n+1}$.
\begin{lemma}\label{grad_v_bound} Consider the update Step 2A in (1.3). Assume $I_H = I_H^1$, we have the following inequality:
\begin{equation}
\begin{gathered}
   \|\nabla v^{n+1}\|^2 \leq \|\nabla \widetilde{v}^{n+1}\|^2 - \|\nabla (v^{n+1}-\widetilde{v}^{n+1})\|^2 \\-k\chi \|I_H^0(\nabla v^{n+1})\|^2
   +k\chi \|I_H^0 (\nabla u(t^{n+1}))\|^2.
   \end{gathered}
\end{equation}
\end{lemma}
\begin{proof}
Take the gradient of Step 2A:
\begin{equation}\label{gradinet_2A}
        \frac{\nabla v^{n+1}-\nabla \widetilde{v}^{n+1}}{k} -\chi I_H^0(\nabla u(t^{n+1})-\nabla v^{n+1})=0.
\end{equation}
Take the inner product of (\ref{gradinet_2A}) with $\nabla v^{n+1}$, apply the polarization identity and multiple by $k$, we have
\begin{equation}
\begin{gathered}
\frac{1}{2}\|\nabla v^{n+1}\|^2 -\frac{1}{2}\|\nabla \widetilde{v}^{n+1}\|^2 + \frac{1}{2}\|\nabla (v^{n+1}-\widetilde{v}^{n+1})\|^2 \\+k\chi \|I_H^0(\nabla v^{n+1})\|^2
 =k\chi (I_H^0(\nabla u(t^{n+1})), I_H^0 (\nabla v^{n+1}))\\
 \leq \frac{k\chi}{2}\|I_H^0(\nabla u(t^{n+1}))\|^2 + \frac{k\chi}{2}\|I_H^0(\nabla v^{n+1})\|^2.
 \end{gathered}
\end{equation}
Multiply the above by 2; we have the final result.
\end{proof}

Define the errors
\[
e^{n}(x)=u(x,t^{n})-v^{n}(x){{\text{ and }}}\widetilde{e}^{n}(x)=u(x,t^{n}%
)-\widetilde{v}^{n}(x).
\]
We use Lemma \ref{lemma_H1_semi_prop} to derive the error estimate in Theorem \ref{error_finite}. 
\begin{lemma}\label{lemma_H1_semi_prop} Consider the update Step 2A in (1.3). Assume $I_H = I_H^1$, and (\ref{I_H_assump1}). Then the following holds:
\begin{equation}\label{grad_e_less_grad_e_tilde}
\|I_H^1e^n\|^2 \geq \frac{1}{2}\|e^n\|^2 - (C_1H)^2\|\nabla e^n\|^2.
\end{equation}
\end{lemma}
\begin{proof}
\begin{equation}
\|I_H^1 e^n\|^2 = \|e^n\|^2 -\|(I-I_H^1) e^n\|^2 - 2(I_H^1e^n, (I-I_H^1)e^n).
\end{equation}
Since
\begin{equation}
\begin{gathered}
(I_H^1e^n, (I-I_H^1)e^n)\leq \|I_H^1e^n\|\|(I-I_H^1)e^n\|\\
\leq \frac{1}{2}\|I_H^1e^n)\|^2 + \frac{1}{2} \|(I-I_H^1)e^n\|^2,
\end{gathered}
\end{equation}
we have
\begin{equation}
   \|I_H^1 e^n\|^2 \geq \frac{1}{2}\|e^n\|^2 -\|(I-I_H^1) e_n\|^2.
\end{equation}
We have the final result by applying (\ref{I_H_assump1}).
\end{proof}

We prove that $\|\nabla e^{n+1}\|\leq \|\nabla \widetilde{e}^{n+1}\|$ next in Lemma \ref{lemma_grad_v}. 
\begin{lemma}\label{lemma_grad_v} Consider Step 2A in (1.3). Assume $I_H = I_H^1$, then the following holds
\begin{equation}
\|\nabla e^{n+1}\|^2 +\|\nabla (e^{n+1}-\widetilde{e}^{n+1})\|^2 + 2k\chi \|I_H^0(\nabla e^{n+1})\|^2 =\|\nabla \widetilde{e}^{n+1}\|^2. 
\end{equation}
\end{lemma}
\begin{proof}
Add and subtract $u^{n+1}$ to Step 2A, and multiple $k$. We have
\begin{equation}\label{s2_e}
e^{n+1}-\widetilde{e}^{n+1} + k\chi I_H^1(e^{n+1})=0.
\end{equation} 
Take the gradient of (\ref{s2_e}). There follows
\begin{equation}
\nabla e^{n+1} - \nabla \widetilde{e}^{n+1} + k \chi \nabla I_H^1(e^{n+1})=0.
\end{equation}
Applying Lemma \ref{layton2002}, yields
\begin{equation}\label{gradient_error}
\nabla e^{n+1} - \nabla \widetilde{e}^{n+1} + k \chi I_H^0(\nabla e^{n+1})=0.
\end{equation}
Take the inner product of (\ref{gradient_error}) with $\nabla e^{n+1}$. Expanding shows that
\begin{equation}
\frac{1}{2} \|\nabla e^{n+1}\|^2 -\frac{1}{2}\|\nabla \widetilde{e}^{n+1}\|^2 +\frac{1}{2}\|\nabla (e^{n+1}-\widetilde{e}^{n+1})\|^2 + k\chi \|I_H^0(\nabla e^{n+1})\|^2 =0.
\end{equation}
Multiplying the above equation by $2$ and rearranging terms, gives (\ref{grad_e_less_grad_e_tilde}).
\end{proof}

We next present the finite time stability in Theorem \ref{stability_v_2A_finite}.

\begin{theorem}\label{stability_v_2A_finite}
Consider Steps 1 and 2A in (1.3). Assume $I_H = I_H^1$, $\nu-6\chi C_1^2H^2>0$ and (\ref{I_H_assump1}). The following holds:
\begin{equation}
\begin{gathered}
\|v^N\|^2  + \sum_{n=0}^{N-1}\|\widetilde{v}^{n+1}-v^{n+1}\|^2+ \sum_{n=0}^{N-1} k \chi \|I_H(v^{n+1})\|^2 \\
+ \sum_{n=0}^{N-1}\frac{1}{2} k \nu \|\nabla \widetilde{v}^{n+1}\|^2 \leq \|v^0\|^2+ \sum_{n=0}^{N-1}\frac{k}{\nu}\|f^{n+1}\|^2_{-1} \\+\sum_{n=0}^{N-1} 3k\chi \|I_H (u(t_{n+1})\|^2
+ \sum_{n=0}^{N-1} 3k^2 \chi^2 C_1^2 H^2 \|I_H^0(\nabla u(t^{n+1}))\|^2.
\end{gathered}
\end{equation}
\end{theorem}
\begin{proof}
Take the inner product of Step 1 with $\widetilde{v}^{n+1}$. This gives,
\begin{equation}\label{s1_inner}
(\frac{\widetilde{v}^{n+1}-v^n}{k}, \widetilde{v}^{n+1}) + \nu (\nabla \widetilde{v}^{n+1}, \nabla \widetilde{v}^{n+1}) + (q^{n+1},\nabla \cdot \widetilde{v}^{n+1}) = (f^{n+1}, \widetilde{v}^{n+1}).
\end{equation}
Multiply (\ref{s1_inner}) by $2k$, apply the polarization identity on $(\widetilde{v}^{n+1}-v^n, \widetilde{v}^{n+1})$, and rearrange terms. We have
\begin{equation}\label{s1_inequ}
\|\widetilde{v}^{n+1}\|^2 -\|v^n\|^2 + \|\widetilde{v}^{n+1}-v^n\|^2 + k\nu \|\nabla \widetilde{v}^{n+1}\|^2 \leq  \frac{k}{\nu}\|f^{n+1}\|^2_{-1}.
\end{equation}
The inner product of Step 2A with $v^{n+1}$ yields
\begin{equation}\label{s2_inner}
(\frac{v^{n+1}-\widetilde{v}^{n+1}}{k}, v^{n+1}) - \chi (I_H(u(t_{n+1})-v^{n+1}), v^{n+1})=0.
\end{equation}
Multiply (\ref{s2_inner}) by $k$ and apply the polarization identity on $(v^{n+1}-\widetilde{v}^{n+1},v^{n+1})$. This gives
\begin{equation}
\begin{gathered}
\frac{1}{2}\|v^{n+1}\|^2 -\frac{1}{2}\|\widetilde{v}^{n+1}\|^2 + \frac{1}{2}\|\widetilde{v}^{n+1}-v^{n+1}\|^2+ k \chi \|I_H(v^{n+1})\|^2 \\
= k \chi (I_H (u(t_{n+1})), I_H(v^{n+1}))-k\chi (I_H(v^{n+1}),(I-I_H)(v^{n+1}))\\
+k\chi (I_H(u(t_{n+1})), (I-I_H)(v^{n+1}))\\
\leq k\chi \|I_H(u(t_{n+1}))\|^2 +\frac{k \chi}{4}\|I_H(v^{n+1})\|^2 + \frac{k \chi}{4}\|I_H(v^{n+1})\|^2 \\
+ k\chi \|(I-I_H)(v^{n+1})\|^2
+ \frac{k\chi}{2}\|I_H (u(t_{n+1}))\|^2 +\frac{k \chi}{2}\|(I-I_H)(v^{n+1})\|^2\\
=\frac{3}{2}k \chi \|I_H(u(t_{n+1}))\|^2 +\frac{k \chi}{2}\|I_H(v^{n+1})\|^2 + \frac{3}{2}k \chi \|(I-I_H)(v^{n+1})\|^2.
\end{gathered}
\end{equation}
Rearranging terms, multiplying the above equation by $2$, and applying (\ref{I_H_assump1}) yields
\begin{equation}
\begin{gathered}
\|v^{n+1}\|^2 -\|\widetilde{v}^{n+1}\|^2 + \|\widetilde{v}^{n+1}-v^{n+1}\|^2+ k \chi \|I_H(v^{n+1})\|^2\\
\leq 3 k \chi \|I_H (u(t_{n+1}))\|^2 + 3k \chi C_1^2 H^2 \|\nabla v^{n+1}\|^2.
\end{gathered}
\end{equation}
By Lemma \ref{grad_v_bound}, we thus have 
\begin{equation}\label{s2_inequ}
\begin{gathered}
\|v^{n+1}\|^2 -\|\widetilde{v}^{n+1}\|^2 + \|\widetilde{v}^{n+1}-v^{n+1}\|^2+ k \chi \|I_H(v^{n+1})\|^2 \\
\leq 3 k \chi \|I_H (u(t_{n+1}))\|^2 + 3k \chi C_1^2 H^2 \left[ \|\nabla \widetilde{v}^{n+1}\|^2 +k\chi \|I_H^0(\nabla u(t^{n+1}))\|^2\right].
\end{gathered}
\end{equation}
Adding (\ref{s1_inequ}) and (\ref{s2_inequ}), we have
\begin{equation}
\begin{gathered}
\|v^{n+1}\|^2 - \|v^n\|^2 + \|\widetilde{v}^{n+1}-v^n\|^2 + \|\widetilde{v}^{n+1}-v^{n+1}\|^2+ k \chi \|I_H(v^{n+1})\|^2  \\+ \frac{1}{2}k \nu \|\nabla \widetilde{v}^{n+1}\|^2+ k (\frac{1}{2}\nu-3\chi C_1^2H^2) \|\nabla \widetilde{v}^{n+1}\|^2 \\
\leq \frac{k}{\nu}\|f^{n+1}\|^2_{-1} + 3k\chi \|I_H (u(t_{n+1})\|^2
+ 3k^2 \chi^2 C_1^2 H^2 \|I_H^0 (\nabla u(t^{n+1}))\|^2.
\end{gathered}
\end{equation}
Since $\nu - 6 C_1 \chi H^2 >0$, we have
\begin{equation}\label{before_sum_stability}
\begin{split}
\|v^{n+1}\|^2 - \|v^n\|^2 + \|\widetilde{v}^{n+1}-v^n\|^2 + \|\widetilde{v}^{n+1}-v^{n+1}\|^2\\+ k \chi \|I_H(v^{n+1})\|^2  + \frac{1}{2}k \nu \|\nabla \widetilde{v}^{n+1}\|^2 \leq \frac{k}{\nu}\|f\|^2_{-1} \\+ 3k\chi \|I_H (u(t_{n+1})\|^2
+ 3k^2 \chi^2 C_1^2 H^2 \|I_H^0(\nabla u(t^{n+1}))\|^2.
\end{split}
\end{equation}
\end{proof}

We next prove infinite time stability using Theorem \ref{stability_v_2A_finite}. We assume that $k \leq \frac{2C_{PF}^2}{\nu}$ (which is typically $\gg 1$) so that we can use the inequality  \begin{equation} 
\ln(1+\frac{k\nu}{2C_{PF}^2}) >\frac{1}{2}(1+ \frac{k\nu}{2C_{PF}^2}) 
\end{equation}
in (\ref{k_tau_condition}) to achieve exponential decay in time corresponding to the initial velocity. 
\begin{theorem}
Consider steps 1 and 2A in (1.3). Assume $I_H = I_H^1$, $\nu-6\chi C_1^2H^2>0$, $k\leq \frac{2C_{PF}^2}{\nu}$, and (\ref{I_H_assump1}). Then, the following holds:
\begin{equation}\label{stability_v_2A_infinite}
\begin{gathered}\|v^{n+1}\|^2\leq \exp\left[-\frac{\nu}{4C_{PF}^2}t_{n+1}\right] \|v^0\|^2 +\frac{2C_{PF}^2}{\nu}\max_{0\leq n'\leq n} F_{n'+1}, \\
\text{ and } \\
 \limsup_{N\to \infty} (\frac{1}{N}\sum_{n=0}^{N-1} \frac{\nu k }{2}\|\nabla v^{n+1}\|^2)
 \leq \\ 
  k \max_{0\leq n\leq N } \left[\frac{1}{\nu}\|f^{n+1}\|^2_{-1} + 3\chi \|I_H (u(t_{n+1})\|^2
+ 3k \chi^2 C_1^2 H^2 \|I_H^0 (\nabla u(t^{n+1}))\|^2 \right].
 \end{gathered}
\end{equation}
\end{theorem}
\begin{proof}
By the Poincar\'e inequality, we have $\|\nabla v\|\geq \frac{1}{C_{PF}}\|v\|.$
Since $\|\nabla v\|\leq \|\nabla \widetilde{v}\|$, we have $\frac{1}{C_{PF}}\|v\|\leq \|\nabla \widetilde{v}\|$. Applying Theorem \ref{stability_v_2A_finite}, we thus have
\begin{equation}\label{before_rearrange_stability}
\begin{gathered}
\left(1+ \frac{k\nu}{2 C_{PF}^2}\right)\|v^{n+1}\|^2 - \|v^n\|^2 \leq \frac{k}{\nu}\|f\|^2_{-1} \\+ 3k\chi \|I_H (u(t_{n+1})\|^2
+ 3k^2 \chi^2 C_1^2 H^2 \|I_H^0(\nabla u(t^{n+1}))\|^2.
\end{gathered}
\end{equation}
We denote $$\tau = \frac{\nu}{2 C_{PF}^2}, \ \  F_{n+1} = \frac{1}{\nu}\|f\|^2_{-1} + 3\chi \|I_H (u(t_{n+1})\|^2
+ 3k \chi^2 C_1^2 H^2 \|I_H^0 (\nabla u(t^{n+1}))\|^2.$$
Thus, we rearrange (\ref{before_rearrange_stability}) to read
\begin{equation}
\begin{gathered}
    \|v^{n+1}\|^2  \leq \frac{1}{1+k\tau} \|v^n\|^2 + \frac{k}{1+k\tau} F_{n+1}\\
   \leq 
    \left(\frac{1}{1+k\tau}\right)^2 \|v^{n-1}\|^2 + \left(\frac{1}{1 +k \tau}\right)^2 k F_n + \frac{k}{1 +k \tau} F_{n+1}.\\
    \end{gathered}
    \end{equation}
By an induction argument, there follows
\begin{equation}
        \begin{gathered}
    \|v^{n+1}\|^2     
   \leq \left(\frac{1}{1+k\tau}\right)^{n+1}\|v^0\|^2 +k \max_{0\leq n'\leq n} F_{n'+1} \sum_{n'=0}^{n}\left( \frac{1}{1+k\tau}\right)^{n'+1}\\
    =(1+k\tau)^{-(n+1)}\|v^0\|^2 +\frac{1}{\tau}\max_{0\leq n'\leq n} F_{n'+1}  \\
  =\exp\left[-(n+1)\ln(1+k\tau)\right]\|v^0\|^2 +\frac{1}{\tau} \max_{0\leq n'\leq n} F_{n'+1}.
    \end{gathered}
\end{equation}
Since we assume $0<k\tau<1$, $\ln(1+k\tau)>\frac{1}{2} k\tau$. Thus, 
\begin{equation}\label{k_tau_condition}
\begin{gathered}
\|v^{n+1}\|^2\leq  \exp\left[-\frac{1}{2}(n+1)k\tau\right] \|v^0\|^2 +\frac{1}{\tau}\max_{0\leq n'\leq n} F_{n'+1}\\
    =\exp\left[-\frac{1}{2}\tau t_{n+1} \right] \|v^0\|^2 +\frac{1}{\tau}\max_{0\leq n'\leq n} F_{n'+1}\\
    =\exp\left[-\frac{\nu}{4C_{PF}^2}t_{n+1}\right] \|v^0\|^2 +\frac{2C_{PF}^2}{\nu}\max_{0\leq n'\leq n} F_{n'+1}.
\end{gathered}
\end{equation}
From Theorem \ref{stability_v_2A_finite}, we have
\begin{equation}
    \|v^{N}\|^2+ \sum_{n=0}^{N-1}\frac{1}{2}k \nu \|\nabla v^{n+1}\|^2 \leq \|v^0\|^2 + k\sum_{n=0}^{N-1} F_{n+1}.
\end{equation}
Divide by $N$, and as $N\to \infty$ we have
\begin{equation}
\begin{gathered}
    \limsup_{N\to \infty} ( \frac{1}{N}\sum_{n=0}^{N-1}\frac{1}{2} k \nu \|\nabla v^{n+1}\|^2) \leq \limsup_{N\to \infty} (\frac{1}{N} \|v^0\|^2 + \frac{1}{N} k \sum_{n=0}^N F_{n+1})\\
    \leq k \max_{0\leq n\leq N} F_{n+1}.
    \end{gathered}
\end{equation}
\end{proof}

We now prove a \textit{first} error estimate for the 2-step nudging (1.3) with 2A below. This is a finite time error estimate (with linear time growth) for the $L^2$ error but an infinite time estimate for time averaged $H^1$ error. 
\begin{theorem}\label{error_finite}
Consider steps 1 and 2A of (1.3). Assume $I_H = I_H^1$, $\nu - 8C_1^2 \chi H^2>0$,  and $\frac{1}{8}\chi -\frac{9C_2^4}{2\nu^3}\|\nabla u^{n+1}\|^4 >0$. We have the following error estimate:
\begin{equation}
\begin{gathered}
\frac{1}{N}\|\widetilde{e}^N\|^2-\frac{1}{N}\|\widetilde{e}^0\|^2 +\frac{1}{N}\sum_{n=0}^{N-1}\left[\|e^n -\widetilde{e}^n\|^2+\|\widetilde{e}^{n+1}-e^n\|^2\right]
\\
+ \frac{\nu k}{2}\frac{1}{N} \sum_{n=0}^{N-1}\|\nabla \widetilde{e}^{n+1}\|^2+ \frac{k\chi}{4} \frac{1}{N}\sum_{n=0}^{N-1} \|e^n\|^2 
\leq \frac{3 k^3 }{\nu}\frac{1}{N} \sum_{n=0}^{N-1} \max_{t_n\leq t\leq t_{n+1}} \|u_{tt}(t)\|^2. 
\end{gathered}
\end{equation}
\end{theorem}
\begin{proof}
We evaluate the NSE at $t_{n+1}$:
\begin{equation}\label{nse}
\frac{u^{n+1}- u^n}{k} + u^{n+1}\cdot \nabla u^{n+1} -\nu \Delta u^{n+1} +\nabla p^{n+1} = f^{n+1} +\rho^{n+1},
\end{equation}
where the consistency error term is $\rho^{n+1}= \frac{u^{n+1}-u^n}{k}-u_t(t_{n+1})$. 

Subtract Step 1 from (\ref{nse}), and then take the inner product with $\widetilde{e}^{n+1}$. We have 
\begin{equation}
\begin{gathered}
(\frac{\widetilde{e}^{n+1}-e^n}{k}, \widetilde{e}^{n+1}) + (u^{n+1} \cdot \nabla u^{n+1}, \widetilde{e}^{n+1})  -(v^n\cdot \nabla  \widetilde{v}^{n+1}, \widetilde{e}^{n+1}) \\+ \nu (\nabla \widetilde{e}^{n+1}, \nabla \widetilde{e}^{n+1}) = (\rho^{n+1}, \widetilde{e}^{n+1}).
\end{gathered}
\end{equation}
Use the polarization identity on $(\frac{\widetilde{e}^{n+1}-e^n}{k}, \widetilde{e}^{n+1})$. We have
\begin{equation}\label{eqn1_error}
\begin{split}
\frac{1}{2k}\left[\|\widetilde{e}^{n+1}\|^2 -\|e^n\|^2+\|\widetilde{e}^{n+1}-e^n\|^2\right] + (u^{n+1}\cdot \nabla u^{n+1}, \widetilde{e}^{n+1})  \\-( v^n\cdot \nabla \widetilde{v}^{n+1}, \widetilde{e}^{n+1})+ \nu (\nabla \widetilde{e}^{n+1}, \nabla \widetilde{e}^{n+1}) = (\rho^{n+1}, \widetilde{e}^{n+1}).
\end{split}
\end{equation}
We calculate $(second$ $eqn, e^{n+1})$, and use the polarization identity to yield 
\begin{equation}
    \frac{1}{2}\left[\|e^{n+1}\|^2 -\|\widetilde{e}^{n+1}\|^2 +\|e^{n+1}-\widetilde{e}^{n+1}\|^2\right] + k \chi (I_H(e^{n+1}), e^{n+1})=0.
\end{equation}
Next, we bound $k \chi (I_H(e^{n+1}), e^{n+1})$.
\begin{equation}
k\chi (I_H(e^{n+1}),e^{n+1}) = k \chi \|I_H(e^{n+1})\|^2 + k\chi (I_H(e^{n+1}), (I-I_H)(e^{n+1})).
\end{equation}
By Cauchy-Schwartz and Young's inequality, there yields
\begin{equation}
\begin{gathered}
k\chi (I_H(e^{n+1}), (I-I_H)(e^{n+1}))\leq \frac{k\chi}{2}\|I_H(e^{n+1})\|^2 + \frac{k \chi}{2}\|(I-I_H)(e^{n+1})\|^2\\
\leq \frac{k\chi}{2}\|I_H(e^{n+1})\|^2 +  \frac{k \chi}{2} C_1^2 H^2 \|\nabla e^{n+1}\|^2.
\end{gathered}
\end{equation}
Multiply by $2$ and rearrange terms. There results
\begin{equation}\label{l2_error_bound}
    \|e^{n+1}\|^2 \leq \|\widetilde{e}^{n+1}\|^2 - \|e^{n+1}-\widetilde{e}^{n+1}\|^2 +  k \chi C_1^2 H^2 \|\nabla e^{n+1}\|^2 -k\chi\|I_H(e^{n+1})\|^2.
\end{equation}
We evaluate (\ref{l2_error_bound}) at $t_{n}$ to replace $\|e^n\|^2$ in  (\ref{eqn1_error}), yielding
\begin{equation}
\begin{gathered}
 \frac{1}{2k} \left[ \|\widetilde{e}^{n+1}\|^2- \|\widetilde{e}^{n}\|^2 + \|\widetilde{e}^{n+1}-e^n\|^2 +\|e^{n}-\widetilde{e}^n\|^2 \right]\\+\frac{1}{2} \chi \|I_H (e^n)\|^2- \frac{\chi}{2} C_1^2 H^2 \|\nabla e^{n}\|^2  + (u^{n+1}\cdot \nabla u^{n+1}, \widetilde{e}^{n+1}) \\-( v^n\cdot \nabla \widetilde{v}^{n+1}, \widetilde{e}^{n+1})\leq (\frac{\widetilde{e}^{n+1}-e^n}{k}, \widetilde{e}^{n+1}).
\end{gathered}
\end{equation}
We now apply Lemma \ref{lemma_H1_semi_prop}, giving
\begin{equation}\label{before_nonlinear_terms}
\begin{gathered}
 \frac{1}{2k} \left[ \|\widetilde{e}^{n+1}\|^2- \|\widetilde{e}^{n}\|^2 + \|\widetilde{e}^{n+1}-e^n\|^2 +\|e^{n}-\widetilde{e}^n\|^2 \right]\\+\frac{\chi}{4}\|e^n\|^2 - \chi C_1^2H^2\|\nabla e^n\|^2  + (u^{n+1}\cdot \nabla u^{n+1}, \widetilde{e}^{n+1})\\
 -( v^n\cdot \nabla \widetilde{v}^{n+1}, \widetilde{e}^{n+1})\leq (\frac{\widetilde{e}^{n+1}-e^n}{k}, \widetilde{e}^{n+1}).
\end{gathered}
\end{equation}
The nonlinear terms are rearranged as
\begin{equation}
\begin{gathered}
(u^{n+1}\cdot \nabla u^{n+1}, \widetilde{e}^{n+1})  -(v^n\cdot \nabla \widetilde{v}^{n+1}, \widetilde{e}^{n+1})
 =(u^{n+1}-u^n\cdot \nabla u^{n+1}, \widetilde{e}^{n+1})\\+ (u^n \cdot \nabla u^{n+1},\widetilde{e}^{n+1})-(v^n \cdot \nabla u^{n+1},\widetilde{e}^{n+1})
+(v^n\cdot \nabla u^{n+1},\widetilde{e}^{n+1})\\-(v^n\cdot \nabla \widetilde{v}^{n+1},\widetilde{e}^{n+1})
=(u^{n+1}-u^n\cdot \nabla u^{n+1}, \widetilde{e}^{n+1}) + (e^n\cdot \nabla  u^{n+1},\widetilde{e}^{n+1}) \\+ (v^n\cdot \nabla  \widetilde{e}^{n+1}, \widetilde{e}^{n+1})
=(u^{n+1}-u^n\cdot \nabla u^{n+1}, \widetilde{e}^{n+1}) +(e^n\cdot \nabla u^{n+1},\widetilde{e}^{n+1}).
\end{gathered}
\end{equation}
We bound the terms $(u^{n+1}-u^n\cdot \nabla u^{n+1}, \widetilde{e}^{n+1})$ and $(e^n\cdot \nabla u^{n+1},\widetilde{e}^{n+1})$ as follows.
\begin{equation}
\begin{gathered}
(u^{n+1}-u^n\cdot \nabla  u^{n+1}, \widetilde{e}^{n+1})\leq C_2\|\nabla (u^{n+1}-u^n)\|\nabla u^{n+1}\|\|\nabla \widetilde{e}^{n+1}\|\\
\leq \frac{\nu}{6} \|\nabla \widetilde{e}^{n+1}\|^2 + \frac{3}{2\nu} C_2^2\|\nabla (u^{n+1}-u^n)\|^2 \|\nabla u^{n+1}\|^2\\
\leq \frac{\nu}{6} \|\nabla \widetilde{e}^{n+1}\|^2 + \frac{3}{2\nu} C_2^2 k \left[\int_{t_n}^{t_{n+1}} \|\nabla u_t\|^2\right]\|\nabla u^{n+1}\|^2\\
 \leq \frac{\nu}{6} \|\nabla \widetilde{e}^{n+1}\|^2 + \frac{3}{2\nu}C_2^2 k^2 \max_{t_n\leq t\leq t_{n+1}} \|\nabla u_t\|^2 \|\nabla u^{n+1}\|^2.
\end{gathered}
\end{equation}
\begin{equation}
\begin{split}
&(e^n\cdot \nabla u^{n+1}, \widetilde{e}^{n+1})\leq C_2 \|e^n\|^{1/2}\|\nabla e^n\|^{1/2} \|\nabla u^{n+1}\|\|\nabla \widetilde{e}^{n+1}\|\\
&\leq \frac{\nu}{6} \|\nabla \widetilde{e}^{n+1}\|^2 +\frac{3C_2^2}{2\nu} \|e^n\|\|\nabla e^{n}\|\|\nabla u^{n+1}\|^2\\
&\leq \frac{\nu}{6} \|\nabla \widetilde{e}^{n+1}\|^2 + \frac{\nu}{8}\|\nabla e^n\|^2 + \frac{9  C_2^4}{2\nu^3}\|\nabla u^{n+1}\|^4 \|e^n\|^2.
\end{split}
\end{equation}
Thus,
\begin{equation}\label{nonlinear_terms}
\begin{gathered}
(u^{n+1}\cdot \nabla u^{n+1}, \widetilde{e}^{n+1})  -(v^n\cdot \nabla \widetilde{v}^{n+1}, \widetilde{e}^{n+1})\leq 
\frac{\nu}{3} \|\nabla \widetilde{e}^{n+1}\|^2 + \frac{\nu}{8}\|\nabla e^n\|^2 \\+ \frac{9 C_2^4}{2\nu^3}\|\nabla u^{n+1}\|^4 \|e^n\|^2+
\frac{3}{2\nu}C_2^2 k^2 \max_{t_n\leq t\leq t_{n+1}} \|\nabla u_t\|^2 \|\nabla u^{n+1}\|^2 .
\end{gathered}
\end{equation}
The consistency error term satisfies
\begin{equation}
\begin{split}
(\rho^{n+1}, \widetilde{e}^{n+1})\leq \|\rho^{n+1}\|\|\nabla \widetilde{e}^{n+1}\|
\leq \frac{3}{2\nu}\|\rho^{n+1}\|^2 +\frac{\nu}{6}\|\nabla \widetilde{e}^{n+1}\|^2.
\end{split}
\end{equation}
By Taylor's Theorem with an integral remainder, we have
\begin{equation}
\rho^{n+1} =\frac{1}{k} \int_{t_n}^{t_{n+1}} u_{tt} (t_n-t') \, dt'.
\end{equation}
Hence,
\begin{equation}
\begin{gathered}
\|\rho^{n+1}\|^2 = \|\frac{1}{k} \int_{t_n}^{t_{n+1}} u_{tt} (t_n-t') \, dt'\|^2
\leq \int_{\Omega}(\int_{t_n}^{t_{n+1}} u_{tt}(t')\, dt')^2\, dx \\
\leq k^2 \max_{t_n\leq t\leq t_{n+1}} \int_{\Omega} u_{tt}^2(t)\, dx
= k^2  \max_{t_n\leq t\leq t_{n+1}} \|u_{tt}(t)\|^2.
\end{gathered}
\end{equation}
Thus,
\begin{equation}\label{consistency_error}
    (\rho^{n+1}, \widetilde{e}^{n+1})\leq \frac{3}{2 \nu}k^2 \max_{t_n\leq t\leq t_{n+1}}\|u_{tt}(t)\|^2 + \frac{\nu}{6}\|\nabla \widetilde{e}^{n+1}\|^2.
\end{equation}
Adding (\ref{before_nonlinear_terms}), (\ref{nonlinear_terms}), and (\ref{consistency_error}), we have
\begin{equation}
\begin{gathered}
\frac{1}{2k} \left[ \|\widetilde{e}^{n+1}\|^2 -\|\widetilde{e}^{n}\|^2 + \|\widetilde{e}^{n+1}-e^n\|^2 +\|e^{n}-\widetilde{e}^n\|^2 \right]\\ 
+ (\frac{\chi}{4}-\frac{9C^4}{2\nu^3}\|\nabla u^{n+1}\|^4 )\|e^n\|^2  + \frac{\nu}{4} \|\nabla \widetilde{e}^{n+1}\|^2 + \frac{\nu}{4} (\|\nabla 
\widetilde{e}^{n+1}\|^2 - \|\nabla e^n\|^2)  \\
+ \frac{1}{8}(\nu -8C_1^2\chi H^2) \|\nabla e^n\|^2
\leq \frac{3}{2\nu}k^2  \max_{t_n\leq t\leq t_{n+1}} \|u_{tt}(t)\|^2.
\end{gathered}
\end{equation}
By Lemma \ref{lemma_grad_v}, $\|\nabla e^n\|\leq \|\nabla \widetilde{e}^n\|$. Hence,
\begin{equation}
\begin{gathered}
\frac{1}{2k} \left[ \|\widetilde{e}^{n+1}\|^2 -\|\widetilde{e}^{n}\|^2 +\|\widetilde{e}^{n+1}-e^n\|^2+\|e^{n}-\widetilde{e}^n\|^2 \right]  \\+ \frac{\nu}{4} \|\nabla \widetilde{e}^{n+1}\|^2   + (\frac{\chi}{4}-\frac{9 C_2^4}{2\nu^3}\|\nabla u^{n+1}\|^4 )\|e^n\|^2  \\
+ \frac{1}{8}(\nu -8C_1^2 \chi H^2) \|\nabla e^n\|^2
\leq \frac{3}{2\nu} k^2  \max_{t_n\leq t\leq t_{n+1}} \|u_{tt}(t)\|^2.
\end{gathered}
\end{equation}
Since we assume that $\nu - 8C_1^2 \chi H^2>0$, and $\frac{1}{8}\chi -\frac{9C_2^4}{2\nu^3}\|\nabla u^{n+1}\|^4 >0$,
\begin{equation}\label{error_before_sum}
\begin{split}
\frac{1}{2k} \left[\|\widetilde{e}^{n+1}\|^2 -\|\widetilde{e}^{n}\|^2 +\|\widetilde{e}^{n+1}-e^n\|^2 +\|e^{n}-\widetilde{e}^n\|^2 \right]\\  + \frac{\nu}{4} \|\nabla \widetilde{e}^{n+1}\|^2 + \frac{\chi}{8}\|e^n\|^2
\leq \frac{3k^2 }{2\nu}  \max_{t_n\leq t\leq t_{n+1}} \|u_{tt}(t)\|^2.
\end{split}
\end{equation}
Multiplying the above inequality by $2k$, and summing over $n$, completes the proof.
\end{proof}

Theorem \ref{error_finite} estimates the error $\widetilde{e}$. We extend it to an analogous estimate of $e$ in Proposition \ref{error_l2_bound}.

\begin{proposition}\label{error_l2_bound} Consider (1.3). Assume $I_H = I_H^1$, $\nu - 8C_1^2 \chi H^2>0$, $\frac{1}{8}\chi -\frac{9C_2^4}{2\nu^3}\|\nabla u^{n+1}\|^4 >0$, and $k\leq \frac{2C_{PF}^2}{\nu}$. We have the following error estimate:
\begin{equation}
\begin{split}
\|e^N\|^2 +\sum_{n=0}^{N-1}\left[\|e^n -\widetilde{e}^n\|^2+\|\widetilde{e}^{n+1}-e^n\|^2\right]
+ \frac{3\nu k}{8} \sum_{n=0}^{N-1}\|\nabla e^{n+1}\|^2\\ + \frac{k\chi}{4} \sum_{n=0}^{N-1} \|e^n\|^2 
\leq \|\widetilde{e}^0\|^2 + \frac{3 k^3 }{\nu} \sum_{n=0}^{N-1} \max_{t_n\leq t\leq t_{n+1}} \|u_{tt}(t)\|^2. 
\end{split}
\end{equation}
\end{proposition}
\begin{proof}
By (\ref{l2_error_bound}), we have
\begin{equation}
    \begin{gathered}
        \|e^{n+1}\|^2 + \frac{3\nu}{8} k\|\nabla e^{n+1}\| \leq \|\widetilde{e}^{n+1}\|^2  +  \frac{3\nu}{8} k\|\nabla \widetilde{e}^{n+1}\|^2 + k \chi C_1^2 H^2 \|\nabla e^{n+1}\|^2\\
        \leq \|\widetilde{e}^{n+1}\|^2 +  \frac{3\nu}{8} k\|\nabla \widetilde{e}^{n+1}\|^2 + k \chi C_1^2 H^2 \|\nabla \widetilde{e}^{n+1}\|^2\\
        \leq \|\widetilde{e}^{n+1}\|^2 +  \frac{3\nu}{8} k\|\nabla \widetilde{e}^{n+1}\|^2 +\frac{\nu k}{8} \|\nabla \widetilde{e}^{n+1}\|^2
        \leq \|\widetilde{e}^{n+1}\|^2 +  \frac{\nu}{2} k\|\nabla \widetilde{e}^{n+1}\|^2.
    \end{gathered}
\end{equation}
Hence, we have the final result.
\end{proof}

Now we prove the infinite time error estimates for $\|e\|$ using Theorem \ref{error_finite}.
\begin{theorem} Consider steps 1 and 2A in (1.3). Assume $I_H = I_H^1$, $\nu - 8C_1^2 \chi H^2>0$, $\frac{1}{8}\chi -\frac{9C_2^4}{2\nu^3}\|\nabla u^{n+1}\|^4 >0$ and $k\leq \frac{2C_{PF}^2}{\nu}$. The following inequalities hold:
\begin{equation}\label{error_l2_infinite}
\frac{4}{5}\|e^{n+1}\|^2 \leq \|\widetilde{e}^{n+1}\|^2 \leq \exp\left[-\frac{\nu}{4C_{PF}^2} t_{n+1} \right]\|\widetilde{e}^0\|^2 + \frac{k^2}{\nu}\frac{2C_{PF}^2}{\nu}, \text{ and }
\end{equation}
\begin{equation}\label{gradient_error_infinite}
\begin{gathered}
  \limsup_{N\to \infty}  \frac{\nu k}{2} \frac{1}{N} \sum_{n=0}^{N-1}\|\nabla e^{n+1}\|^2 \leq   \limsup_{N\to \infty}  \frac{\nu k}{2} \frac{1}{N} \sum_{n=0}^{N-1}\|\nabla \widetilde{e}^{n+1}\|^2 
  \\
  \leq \frac{3k^3}{\nu} \max_{0 \leq t \leq Nk} \|u_{tt}(t)\|^2. 
  \end{gathered}
\end{equation}
\end{theorem}

\begin{proof}
By (\ref{error_before_sum}), we have
\begin{equation}\label{error_inequality}
\|\widetilde{e}^{n+1}\|^2 -\|\widetilde{e}^n\|^2 +\frac{k\nu}{2}\|\nabla \widetilde{e}^{n+1}\|^2 \leq \frac{3k^3}{\nu}\max_{t_n \leq t\leq t_{n+1}} \|u_{tt}(t)\|^2.
\end{equation}
By the Poincar\'e inequality, we have
\begin{equation}
    \|\nabla \widetilde{e}\|\geq \frac{1}{C_{PF}}\|\widetilde{e}\|.
\end{equation}
Thus, (\ref{error_inequality}) becomes
\begin{equation}\label{eqn_4_52}
    (1+\frac{k\nu}{2C_{PF}^2})\|\widetilde{e}^{n+1}\|\leq \|\widetilde{e}^n\|^2 + \frac{3k^3}{\nu}\max_{t_n \leq t \leq t_{n+1}} \|u_{tt}(t)\|^2.
\end{equation}
Denote $\tau = \frac{\nu}{2C_{PF}^2}$, divide $1+k\tau$ to both sides of (\ref{eqn_4_52}). We have
\begin{equation}
\begin{gathered}
\|\widetilde{e}^{n+1}\|^2 \leq \frac{1}{1+k \tau}\|\widetilde{e}^n\|^2 + \frac{k}{1+k\tau}\frac{k^2}{\nu} \max_{t_n\leq t_{n+1}} \|u_{tt}(t)\|^2\\
\leq (\frac{1}{1+k\tau})^2 \|\widetilde{e}^{n-1}\|^2 +  (\frac{1}{1+k\tau})^2\frac{k^2}{\nu}\max_{t_{n-1}\leq t\leq t_n}\|u_{tt}\|^2
\\
+ \frac{k}{1+k\tau}\frac{k^2}{\nu} \max_{t_n \leq tt_{n+1}} \|u_{tt}(t)\|^2.\\
\end{gathered}
\end{equation}
By a simple induction argument, there follows
\begin{equation}
    \begin{gathered}
\|\widetilde{e}^{n+1}\|^2 \leq (\frac{1}{1+k\tau})^{n+1} \|\widetilde{e}^{0}\|^2+ \frac{k^2}{\nu}\max_{0\leq t\leq t_{n+1}}\|u_{tt}\|^2 k \sum_{n'=0}^{n} (\frac{1}{1+k\tau})^{n'+1}\\
 = \exp\left[ -(n+1)\ln (k\tau)\right]\|\widetilde{e}^0\|^2 + \frac{k^2}{\nu}\frac{1}{\tau}.
\end{gathered}
\end{equation}
Since $0<k\tau<1$, $\ln(1+k\tau)>\frac{1}{2} k\tau$, we thus have
\begin{equation}
\|\widetilde{e}^{n+1}\|^2 \leq \exp\left[-\frac{\tau}{2} t_{n+1}\right]\|\widetilde{e}^0\|^2 + \frac{k^2}{\nu}\frac{1}{\tau}.
\end{equation}
By (\ref{l2_error_bound}), we have
\begin{equation}
    \begin{gathered}
        \|e^{n+1}\|^2\leq \|\widetilde{e}^{n+1}\|^2  + k \chi C_1^2 H^2 \|\nabla e^{n+1}\|^2
        \leq \|\widetilde{e}^{n+1}\|^2  + k \chi C_1^2 H^2 \|\nabla \widetilde{e}^{n+1}\|^2\\
        \leq \|\widetilde{e}^{n+1}\|^2 +\frac{\nu k}{8} \|\nabla \widetilde{e}^{n+1}\|^2
        \leq \|\widetilde{e}^{n+1}\|^2 + \frac{\nu k}{8C_{PF}^2}\|\widetilde{e}^{n+1}\|^2 
        \leq \frac{5}{4}\|\widetilde{e}^{n+1}\|^2.
    \end{gathered}
\end{equation}
Thus we have (\ref{error_l2_infinite}). By Theorem \ref{error_finite},
\begin{equation}
\begin{split}
\|\widetilde{e}^N\|^2
+ \frac{\nu k}{2} \sum_{n=0}^{N-1}\|\nabla \widetilde{e}^{n+1}\|^2
\leq \|\widetilde{e}^0\|^2+\frac{3 k^3 }{\nu} \sum_{n=0}^{N-1} \max_{t_n\leq t\leq t_{n+1}} \|u_{tt}(t)\|^2. 
\end{split}
\end{equation}
Divide by $N$, and let $N\to \infty$. We have
\begin{equation}
\begin{gathered}
  \limsup_{N\to \infty} \frac{1}{N} \frac{\nu k}{2} \sum_{n=0}^{N-1}\|\nabla \widetilde{e}^{n+1}\|^2 
  \\
  \leq \limsup_{N\to \infty} \left[\frac{1}{N}\|\widetilde{e}^0\|^2 + \frac{3k^3}{\nu}\frac{1}{N} \sum_{n=0}^{N-1} \max_{t_n\leq t\leq t_{n+1}} \|u_{tt}(t)\|^2\right] \\
  \leq \frac{3k^3}{\nu} \max_{0 \leq n \leq N} \|u_{tt}(t)\|^2. 
  \end{gathered}
\end{equation}
Since $\|\nabla e^{n+1}\|\leq \|\nabla \widetilde{e}^{n+1}\|$, we have (\ref
{gradient_error_infinite}).
\end{proof}

\begin{remark}
The following method of estimating error may generalize to the 2-step method. For problems with smooth solutions (e.g., test problems constructed by the method of manufactured solutions) the $H-$condition can be altered as
follows. For a vector function $w(x)$ define $\lambda _{T}(w)$, the microscale of $w$, as 
\begin{equation*}
\lambda _{T}(w):=\frac{\|w\|}{\|\nabla w\|},d=2{{{{\text { or }}}}}3.
\end{equation*}
Here $\lambda _{T}$\ represents an average length-scale,
analogous to the Taylor micro-scale (which would also involve time averaging). Then $\|\nabla e\|\leq \lambda _{T}(e)^{-1}\|e\|$. The term that
gives rise to the $H$-condition is now estimated by 
\begin{equation*}
\chi \|(I-I_{H})e\|^{2}\leq \chi \left (C_{1}H\|\nabla e\|\right )^{2}\leq \chi C_{1}^{2}H^{2}\lambda _{T}^{-2}\|e\|^{2}.
\end{equation*}
It can now be subsumed into the $\chi $-condition (in 3d here). This gives long time estimates when $H$ is small with respect to $\lambda _{T}(e)$ (rather than ${{{{\mathcal{R}}}}}e$) from the inequality. 
\begin{gather*}
\frac{d}{dt}\|e\|^{2}+\nu \|\nabla e\|^{2} \\
+2\left [\chi \left (1-C_{1}^{2}\left (\frac{H}{\lambda _{T}(e)}\right
)^{2}\right )-\frac{2048}{19683}\nu ^{-3}\|\nabla u\|^{4}\right
]\|e\|^{2}\leq 0.
\end{gather*}
\end{remark}

\begin{remark}
When $I_H = I_H^0$ in Step 2A, the results of this section are open problems. The key is to prove an estimate $\|\nabla e^{n+1}\|\leq C \|\nabla \widetilde{e}^{n+1}\|$. At this point we can only prove this for timestep $k$ small. We have then 
\begin{equation}
\begin{split}
    \|\nabla e\|\leq \|\nabla \widetilde{e}\|+ k\chi \|\nabla I_He\|
    \leq \|\nabla \widetilde{e}\| + k \chi C \|\nabla e\|,
\end{split}
\end{equation}
yielding the required bound for $k \chi C <1$.
\end{remark}

\section{Analysis of Step 2B}

We present the infinite time error analysis for Step 2B when $\nu _{turb}=0$ and $I_H$ is an $L^2$ projection, a case not covered in Section 4. 
 
\begin{gather}
\begin{array}{ccc}
Step{{{\text { }}}}1 &  & 
\begin{array}{c}
\frac{\widetilde {v}^{n+1}-v^{n}}{k}+v^{n}\cdot \nabla \widetilde {v}%
^{n+1}-\nu \triangle \widetilde {v}^{n+1}+\nabla q^{n+1}=f(x){{{{\text { }}}}%
} \\ 
\nabla \cdot \widetilde {v}^{n+1}=0,{{{{\text { in }}}}}\Omega ,%
\end{array}
\\ 
Step & 2B & \frac{v^{n+1}-\widetilde {v}^{n+1}}{k}-\nu \triangle \left
(v^{n+1}-\widetilde {v}^{n+1}\right )-\chi I_{H}(u(t^{n+1})-v^{n+1})=0,%
\end{array}
\\
{{{\text {where }}}}I_{H}u\in X^{H}{{{\text { satisfies: }}}}%
(u-I_{H}u,w^{H})=0{{{\text { , for all }}}}w^{H}\in X^{H}.  \notag
\end{gather}

Assume that $I_{H}$\ is an $L^{2}$ projection satisfying 
\begin{equation}
\|(I-I_{H})w\|\leq C_{1}H\|\nabla w\|{{{{\text { for all }}}}}w\in \left
(H_{0}^{1}(\Omega )\right )^{d},d=2{{{{\text { or }}}}}3.
\end{equation}
Suppose $\chi $ is large enough and $H$ is small enough that 
\begin{equation*}
\chi -\frac{2}{\nu }\|\nabla u(t^{n+1},x)\|_{\infty }^{2}>0{{\text { }}}and{{%
\text { }}}\nu -2\chi \left (C_{1}H\right )^{2}>0.
\end{equation*}
Then Step 1 with Step 2B has an infinite predictability horizon.

Set $e:=u-v,\widetilde{e}:=u-\widetilde{v}$. Subtraction gives the error
equations 
\begin{gather*}
\frac{\widetilde{e}^{n+1}-e^{n}}{k}+\left[ u^{n}\cdot \nabla
u^{n+1}-v^{n}\cdot \nabla \widetilde{v}^{n+1}\right] -\nu \triangle 
\widetilde{e}^{n+1}+\nabla (p-q)^{n+1}=\rho ^{n+1} \\
\rho ^{n+1}={{{\text{consistency error term}}}}=\frac{u(t^{n+1},x)-u(t^{n},x)%
}{k}-u_{t}(t^{n+1},x), \\
e^{n+1}-\widetilde{e}^{n+1}-k\nu \triangle (e_{n+1}-\widetilde{e}%
_{n+1})+k\chi I_{H}\left( e^{n+1}\right) =0.
\end{gather*}%
We calculate $<first$ $eqn,\widetilde{e}^{n+1}>$ and $<second$ $eqn,e^{n+1}>$, use the polarization identity and integrate by parts, several times. This yields 
\begin{gather*}
\frac{1}{2}\left[ \|\widetilde{e}^{n+1}\|^{2}-\|e^{n}\|^{2}+\|\widetilde{e}%
^{n+1}-e^{n}\|^{2}\right] + \\
k(u^{n}\cdot \nabla u^{n+1}-v^{n}\cdot \nabla \widetilde{v}^{n+1},\widetilde{%
e}^{n+1})+k\nu \|\nabla \widetilde{e}^{n+1}\|^{2}=k(\rho ^{n+1},\widetilde{e}%
^{n+1}). \\
{{{\text{Rearranging:}}}} \\
\frac{1}{2}\left[ \|e_{n+1}\|^{2}+k\nu \|\nabla e_{n+1}\|^{2}\right] -\frac{1%
}{2}\left[ \|\widetilde{e}_{n+1}^{h}\|^{2}+k\nu \|\nabla \widetilde{e}%
_{n+1}^{h}\|^{2}\right] + \\
+\frac{1}{2}\left[ \|e_{n+1}^{h}-\widetilde{e}_{n+1}^{h}\|^{2}+k\nu \|\nabla
(e_{n+1}^{h}-\widetilde{e}_{n+1}^{h})\|^{2}\right] +k\chi
\|I_{H}e^{n+1}\|^{2}=0.
\end{gather*}%
The nonlinear terms satisfy 
\begin{gather*}
(u^{n}\cdot \nabla u^{n+1}-v^{n}\cdot \nabla \widetilde{v}^{n+1},\widetilde{e%
}^{n+1})= \\
=(u^{n}\cdot \nabla u^{n+1}-v^{n}\cdot \nabla u^{n+1}+v^{n}\cdot \nabla
u^{n+1}-\widetilde{v}^{n+1}\cdot \nabla \widetilde{v}^{n+1},\widetilde{e}%
^{n+1}) \\
=(e^{n}\cdot \nabla u^{n+1}-v^{n}\cdot \nabla \widetilde{e}^{n+1},\widetilde{%
e}^{n+1}) \\
=(e^{n}\cdot \nabla u^{n+1},\widetilde{e}^{n+1})\leq \frac{\nu }{4}\|\nabla 
\widetilde{e}^{n+1}\|^{2}+\frac{1}{\nu }\|\nabla u^{n+1}\|_{\infty
}^{2}\|e^{n}\|^{2}.
\end{gather*}%
Thus, Steps 1 and 2 imply, respectively, 
\begin{gather*}
\frac{1}{2}\left[ \|\widetilde{e}^{n+1}\|^{2}-\|e^{n}\|^{2}+\|\widetilde{e}%
^{n+1}-e^{n}\|^{2}\right] \\
+k\frac{3\nu }{4}\|\nabla \widetilde{e}^{n+1}\|^{2}-k\frac{1}{\nu }\|\nabla
u^{n+1}\|_{\infty }^{2}\|e^{n}\|^{2}\leq k(\tau ^{n+1},\widetilde{e}^{n+1}),{%
{{\text{ and}}}} \\
\frac{1}{2}\left[ \|e_{n+1}\|^{2}+k\nu \|\nabla e_{n+1}\|^{2}\right] -\frac{1%
}{2}\left[ \|\widetilde{e}_{n+1}^{h}\|^{2}+k\nu \|\nabla \widetilde{e}%
_{n+1}^{h}\|^{2}\right] \\
+\frac{1}{2}\left[ \|e_{n+1}^{h}-\widetilde{e}_{n+1}^{h}\|^{2}+k\nu \|\nabla
(e_{n+1}^{h}-\widetilde{e}_{n+1}^{h})\|^{2}\right] +k\chi
\|I_{H}e^{n+1}\|^{2}=0.
\end{gather*}%
Adding gives 
\begin{gather*}
\frac{1}{2}\left[ \|e^{n+1}\|^{2}\right] -\frac{1}{2}\left[ \|e^{n}\|^{2}%
\right] \\
\frac{1}{2}\left[ \|\widetilde{e}^{n+1}-e^{n}\|^{2}\right] +\frac{k\nu }{2}%
\|\nabla e_{n+1}\|^{2}+k\frac{\nu }{4}\|\nabla \widetilde{e}^{n+1}\|^{2} \\
\frac{1}{2}\left[ \|e_{n+1}^{h}-\widetilde{e}_{n+1}^{h}\|^{2}+k\nu \|\nabla
(e_{n+1}^{h}-\widetilde{e}_{n+1}^{h})\|^{2}\right] +k\chi
\|I_{H}e^{n+1}\|^{2} \\
-k\frac{1}{\nu }\|\nabla u^{n+1}\|_{\infty }^{2}\|e^{n}\|^{2}\leq k(\rho
^{n+1},\widetilde{e}^{n+1}).
\end{gather*}%
Consider the term $+k\chi $ $\|I_{H}e^{n+1}$ $\|^{2}$. Rewrite it as 
\begin{align*}
k\chi \|I_{H}e^{n+1}\|^{2}& =\frac{1}{2}k\chi \|I_{H}e^{n+1}\|^{2}+\frac{1}{2%
}k\chi \|I_{H}e^{n}\|^{2}+\frac{1}{2}k\chi \|I_{H}e^{n+1}\|^{2} \\
& -\frac{1}{2}k\chi \|I_{H}e^{n}\|^{2}.
\end{align*}%
Thus 
\begin{gather*}
\frac{1}{2}\left[ \|e^{n+1}\|^{2}+\frac{1}{2}k\chi \|I_{H}e^{n+1}\|^{2}%
\right] -\frac{1}{2}\left[ \|e^{n}\|^{2}+\frac{1}{2}k\chi \|I_{H}e^{n}\|^{2}%
\right] \\
\frac{1}{2}\left[ \|\widetilde{e}^{n+1}-e^{n}\|^{2}\right] +\frac{k\nu }{2}%
\|\nabla e_{n+1}\|^{2}+k\frac{\nu }{4}\|\nabla \widetilde{e}^{n+1}\|^{2} \\
\frac{1}{2}\left[ \|e_{n+1}^{h}-\widetilde{e}_{n+1}^{h}\|^{2}+k\nu \|\nabla
(e_{n+1}^{h}-\widetilde{e}_{n+1}^{h})\|^{2}\right] + \\
+\frac{k\chi }{2}\left[ \|I_{H}e^{n+1}\|^{2}+\|I_{H}e^{n}\|^{2}\right] -k%
\frac{1}{\nu }\|\nabla u^{n+1}\|_{\infty }^{2}\|e^{n}\|^{2}\leq k(\rho
^{n+1},\widetilde{e}^{n+1}).
\end{gather*}%
Note that 
\begin{gather*}
\|e^{n}\|^{2}=\|I_{H}e^{n}\|^{2}+\|(I-I_{H})e^{n}\|^{2}\leq
\|I_{H}e^{n}\|^{2}+\left( C_{1}H\|\nabla e^{n}\|\right) ^{2},thus \\
\frac{1}{2}k\chi \|I_{H}e^{n}\|^{2}\geq \frac{1}{2}k\chi \|e^{n}\|^{2}-\frac{%
1}{2}k\chi \left( C_{1}H\right) ^{2}\|\nabla e^{n}\|^{2},
\end{gather*}%
so that we have 
\begin{gather*}
\frac{1}{2}\left[ \|e^{n+1}\|^{2}+\frac{1}{2}k\chi \|I_{H}e^{n+1}\|^{2}%
\right] -\frac{1}{2}\left[ \|e^{n}\|^{2}+\frac{1}{2}k\chi \|I_{H}e^{n}\|^{2}%
\right] \\
\frac{1}{2}\left[ \|\widetilde{e}^{n+1}-e^{n}\|^{2}\right] +\frac{k\nu }{2}%
\|\nabla e_{n+1}\|^{2}+k\frac{\nu }{4}\|\nabla \widetilde{e}^{n+1}\|^{2} \\
\frac{1}{2}\left[ \|e_{n+1}^{h}-\widetilde{e}_{n+1}^{h}\|^{2}+k\nu \|\nabla
(e_{n+1}^{h}-\widetilde{e}_{n+1}^{h})\|^{2}\right] + \\
+\frac{k\chi }{2}\|I_{H}e^{n+1}\|^{2}+\frac{1}{2}k\chi \|e^{n}\|^{2}-\frac{1%
}{2}k\chi \left( CH\right) ^{2}\|\nabla e^{n}\|^{2} \\
-k\frac{1}{\nu }\|\nabla u^{n+1}\|_{\infty }^{2}\|e^{n}\|^{2}\leq k(\rho
^{n+1},\widetilde{e}^{n+1}).
\end{gather*}%
Regrouping gives
\begin{gather*}
\frac{1}{2}\left[ \|e^{n+1}\|^{2}+\frac{1}{2}k\chi \|I_{H}e^{n+1}\|^{2}%
\right] -\frac{1}{2}\left[ \|e^{n}\|^{2}+\frac{1}{2}k\chi \|I_{H}e^{n}\|^{2}%
\right] \\
\frac{1}{2}\left[ \|\widetilde{e}^{n+1}-e^{n}\|^{2}\right] +\frac{k\nu }{2}%
\|\nabla e_{n+1}\|^{2}+k\frac{\nu }{4}\|\nabla \widetilde{e}^{n+1}\|^{2} \\
\frac{1}{2}\left[ \|e_{n+1}^{h}-\widetilde{e}_{n+1}^{h}\|^{2}+k\nu \|\nabla
(e_{n+1}^{h}-\widetilde{e}_{n+1}^{h})\|^{2}\right] + \\
+\frac{k\chi }{2}\|I_{H}e^{n+1}\|^{2}+k\left[ \frac{\chi }{2}-\frac{1}{\nu }%
\|\nabla u^{n+1}\|_{\infty }^{2}\right] \|e^{n}\|^{2} \\
-\frac{1}{2}k\chi \left( C_{1}H\right) ^{2}\|\nabla e^{n}\|^{2}\leq k(\rho
^{n+1},\widetilde{e}^{n+1}).
\end{gather*}%
There is one remaining negative term: $-\frac{1}{2}k\chi \left(
C_{1}H\right) ^{2}\|\nabla e^{n}\|^{2}$ . We proceed as above. Write 
\begin{equation*}
\frac{k\nu }{2}\|\nabla e_{n+1}\|^{2}=\frac{k\nu }{4}\|\nabla e_{n+1}\|^{2}+%
\frac{k\nu }{4}\|\nabla e_{n}\|^{2}+\frac{k\nu }{4}\|\nabla e_{n+1}\|^{2}-%
\frac{k\nu }{4}\|\nabla e_{n}\|^{2}.
\end{equation*}%
Use this and regroup. This gives 
\begin{gather*}
\frac{1}{2}\left[ \|e^{n+1}\|^{2}+\frac{1}{2}k\chi \|I_{H}e^{n+1}\|^{2}+%
\frac{k\nu }{4}\|\nabla e_{n+1}\|^{2}\right] \\
-\frac{1}{2}\left[ \|e^{n}\|^{2}+\frac{1}{2}k\chi \|I_{H}e^{n}\|^{2}+\frac{%
k\nu }{4}\|\nabla e_{n}\|^{2}\right] \\
\frac{1}{2}\left[ \|\widetilde{e}^{n+1}-e^{n}\|^{2}\right] +\frac{k}{4}\left[
\nu \|\nabla e_{n+1}\|^{2}+\nu \|\nabla \widetilde{e}^{n+1}\|^{2}\right] + \\
\frac{1}{2}\left[ \|e_{n+1}^{h}-\widetilde{e}_{n+1}^{h}\|^{2}+k\nu \|\nabla
(e_{n+1}^{h}-\widetilde{e}_{n+1}^{h})\|^{2}\right] \\
+\frac{k\chi }{2}\|I_{H}e^{n+1}\|^{2}+k\left[ \frac{\chi }{2}-\frac{1}{\nu }%
\|\nabla u^{n+1}\|_{\infty }^{2}\right] \|e^{n}\|^{2} \\
+\frac{k}{2}\left( \frac{\nu }{2}-\chi \left( CH\right) ^{2}\right) \|\nabla
e^{n}\|^{2}\leq k(\rho ^{n+1},\widetilde{e}^{n+1}).
\end{gather*}%
The proof is then completed with standard arguments in the numerical analysis of the NSE. This yields a uniform in time error estimate (and thus an infinite predictability horizon).

\section{Computational Explorations}
In this section, we present numerical experiments to validate the performance of the 2-step nudging method. Section~\ref{rates_of_conv} verifies the theoretical temporal convergence rates using a test case with a known exact solution. In Section~\ref{test:offset_cyl}, we explore a more complex 2d flow scenario and evaluate the effectiveness of the 2-step nudging by comparing it with standard nudging and a baseline simulation without data assimilation.

\subsection{Convergence rates in time}\label{rates_of_conv}
We verify the theoretical temporal convergence of the backward Euler time discretization of the 2-step nudging using an exact solution from~\cite{CFLS25} on the domain \( \Omega = (0,1) \times (0,1) \). The exact velocity and pressure are given by
\[
u(x,y,t) = \exp(t)(\cos y, \sin x)^\top, \quad p(x,y,t) = (x - y)(1 + t).
\]
Note that $\nabla \times (u\cdot \nabla u)=0$ so, like most exact NSE solutions, this will be a solution of the Stokes problem with modified pressure. These are inserted in the NSE to calculate the corresponding body force $f(x,t)$ with viscosity $\nu=1$. To isolate temporal errors, we use a fine mesh with $43266$ degrees of freedom (dof), minimizing spatial error effects. The final time is $T=2$. Here we choose $I_H$ as the $L^2$ projection. We use a Scott-Vogelius finite element pair with a barycenter refined mesh. We choose the time-step sizes $k = 1/4, 1/8, 1/16, 1/32$, and $1/64$. The $L^2$ velocity errors at $T=2$ for the 2-step nudging and standard nudging are presented in Table \ref{tab:error-main}. We observed a first-order convergence rate for both methods. 
\begin{table}[H]
    \centering
    \begin{subtable}{0.3\linewidth} %
     \centering
        \begin{tabular}{|c|c|c|}
            \hline
            $k$ & $\|u-v\|$ & rate  \\
            \hline
            1/4       & 0.00214  & -      \\
            1/8     & 0.00121   & 0.82   \\
            1/16     & 0.00065  & 0.90\\
            1/32     & 0.00034   & 0.95    \\
            1/64    & 0.000170   & 0.98   \\
            \hline
        \end{tabular}
        \caption{2-step nudging (2A).}
        \label{tab:accu2A}
    \end{subtable}
    \hfill
        \begin{subtable}{0.3\linewidth} %
     \centering
        \begin{tabular}{|c|c|c|}
            \hline
            $k$ & $\|u-v\|$ & rate  \\
            \hline
            1/4       & 0.00265  & -      \\
            1/8     &  0.00136  &  0.96  \\
            1/16     &  0.00069 &   0.98 \\
            1/32     &  0.00034 &   1.02  \\
            1/64    &   0.000173&   0.97 \\
            \hline
        \end{tabular}
        \caption{2-step nudging (2B).}
        \label{tab:accu2B}
    \end{subtable}
    \hfill
    \begin{subtable}{0.3\linewidth}
    \centering
 \begin{tabular}{|c|c|c|}
            \hline
            $k$ & $\|u-v\|$ & rate \\
            \hline
            1/4       & 0.00266   & -      \\
            1/8     & 0.00140   & 0.97   \\
            1/16     & 0.00069  & 0.98   \\
            1/32     & 0.00035  & 0.99   \\
            1/64    & 0.000173   & 1.0   \\
            \hline
        \end{tabular}
        \caption{Standard nudging.}
        \label{tab:accu1}
    \end{subtable}
    \caption{Velocity errors (in $L^2$ norm) and convergence rates for 2-step nudging with 2A and 2B, and standard nudging. All methods exhibit first-order accuracy in time.}
    \label{tab:error-main}
\end{table}
\subsection{Flow between offset cylinders}\label{test:offset_cyl}
\begin{figure}[h]
    \centering
    \includegraphics[width=0.9\linewidth]{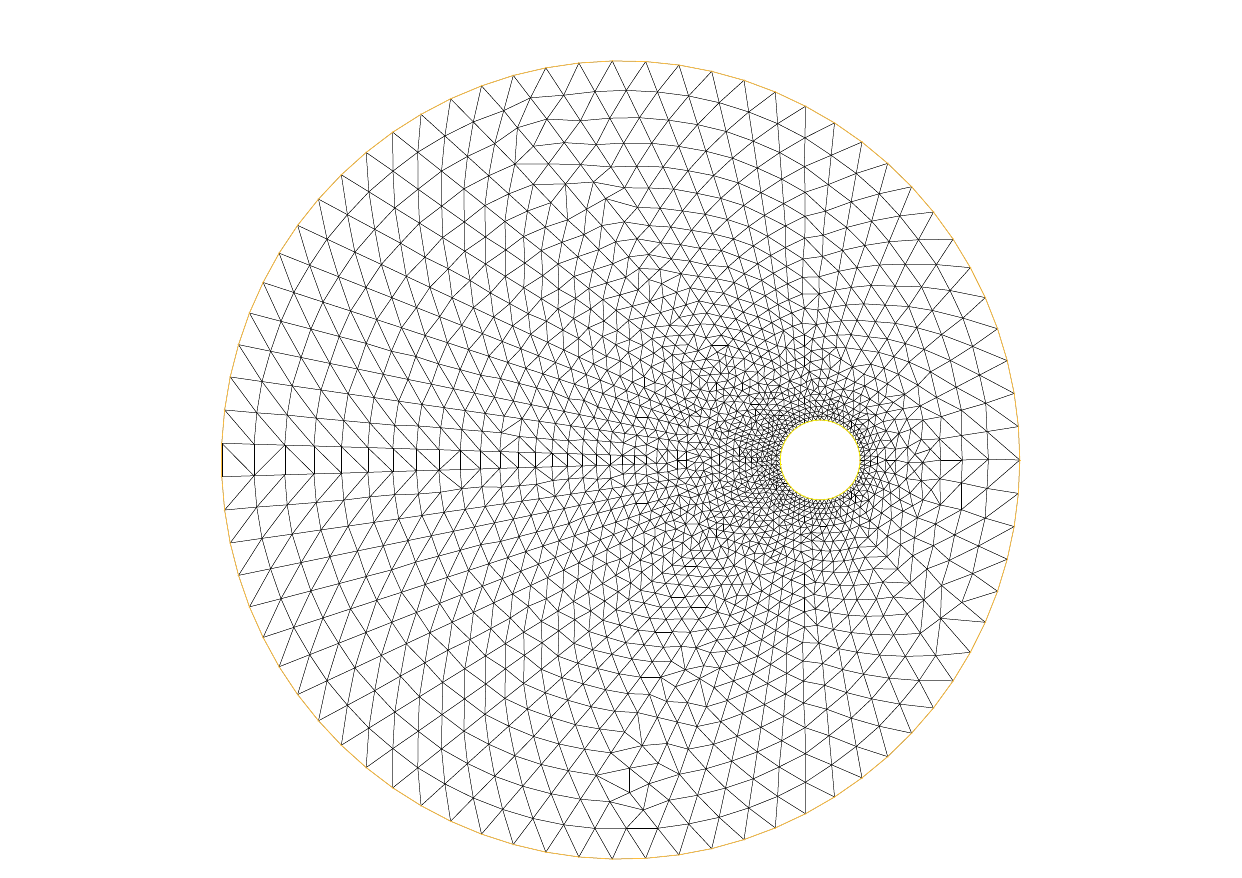}
    \caption{Computational domain with a coarse mesh, consisting of 75 points along the outer boundary and 60 points along the inner boundary.}
    \label{fig:offset_coarser_mesh}
\end{figure}

We evaluate the 2-step nudging methods' performance by using a complex two-dimensional flow at a higher Reynolds number. The computational domain is a disk containing a small, off-center circular obstacle. Let \( r_1 \) and \( r_2 \) denote the radii of the outer and inner circles, respectively, with $r_1 = 1$ and $r_2 = 0.1$. The inner circle (the obstacle) is centered at \( \mathbf{c} = (c_1, c_2) = \left(\frac{1}{2}, 0\right) \). The domain is defined as
\[
\Omega = \left\{ (x, y) \in \mathbb{R}^2 : x^2 + y^2 \leq r_1^2 \text{ and } (x - c_1)^2 + (y - c_2)^2 \geq r_2^2 \right\}.
\]
Figure \ref{fig:offset_coarser_mesh} presents the computational domain with a coarse mesh, generated using Delaunay triangulation. The mesh consists of $75$ points along the outer boundary and $60$ points along the inner boundary. The flow is driven by a counterclockwise rational body force $f$:
\begin{equation}
    f(x,y,t) = \min(1,t) (- 4y(1-x^2-y^2), 4x(1-x^2-y^2))^\top.
\end{equation}
We impose no-slip Dirichlet boundary conditions on the velocity field. The final simulation time is set to $T = 25$, with a time step size $k = 0.01$. The kinematic viscosity is $\nu = 10^{-3}$. We set the characteristic length $L = 1$ and velocity $U = 1$, which gives a Reynolds number of $\mathrm{Re} = \frac{UL}{\nu} = 1000$. The initial condition for the true velocity is $u_0 \equiv 0$ and for the nudged system (with $\mathcal{O}(k)$ error) is
\[
v_0 = k \left(0.1^2- (x- 0.5)^2 - y^2 \right) \left( 1-x^2-y^2\right) \left(-4y,4x\right)^\top.
\]
This $v_0$ is non-zero and consistent with the prescribed boundary conditions. We use backward Euler for the time discretization and Taylor-Hood ($P2-P1$) for the velocity and pressure. The nonlinear convection term is treated using a skew-symmetric formulation with divergence correction:
\begin{equation}
    (v^n\cdot \nabla \widetilde{v}^{n+1}, w) +\frac{1}{2} (\nabla \cdot v^n) (\widetilde{v}^{n+1}, w), \forall w \in X^H.
\end{equation}
We also compare against the standard nudging method, using the same settings as in the two-step nudging case.

To approximate the true solution \( u \), we perform a direct numerical simulation (DNS) using a finer mesh, consisting of 120 points along the outer boundary and 96 points along the inner boundary. For this DNS reference solution, we use the BDF2 (second-order Backward Differentiation Formula) time-stepping scheme. The nonlinear term is also treated using the skew-symmetric formulation with divergence correction.

\begin{figure}[h]
    \centering
    \includegraphics[width=0.9\linewidth]{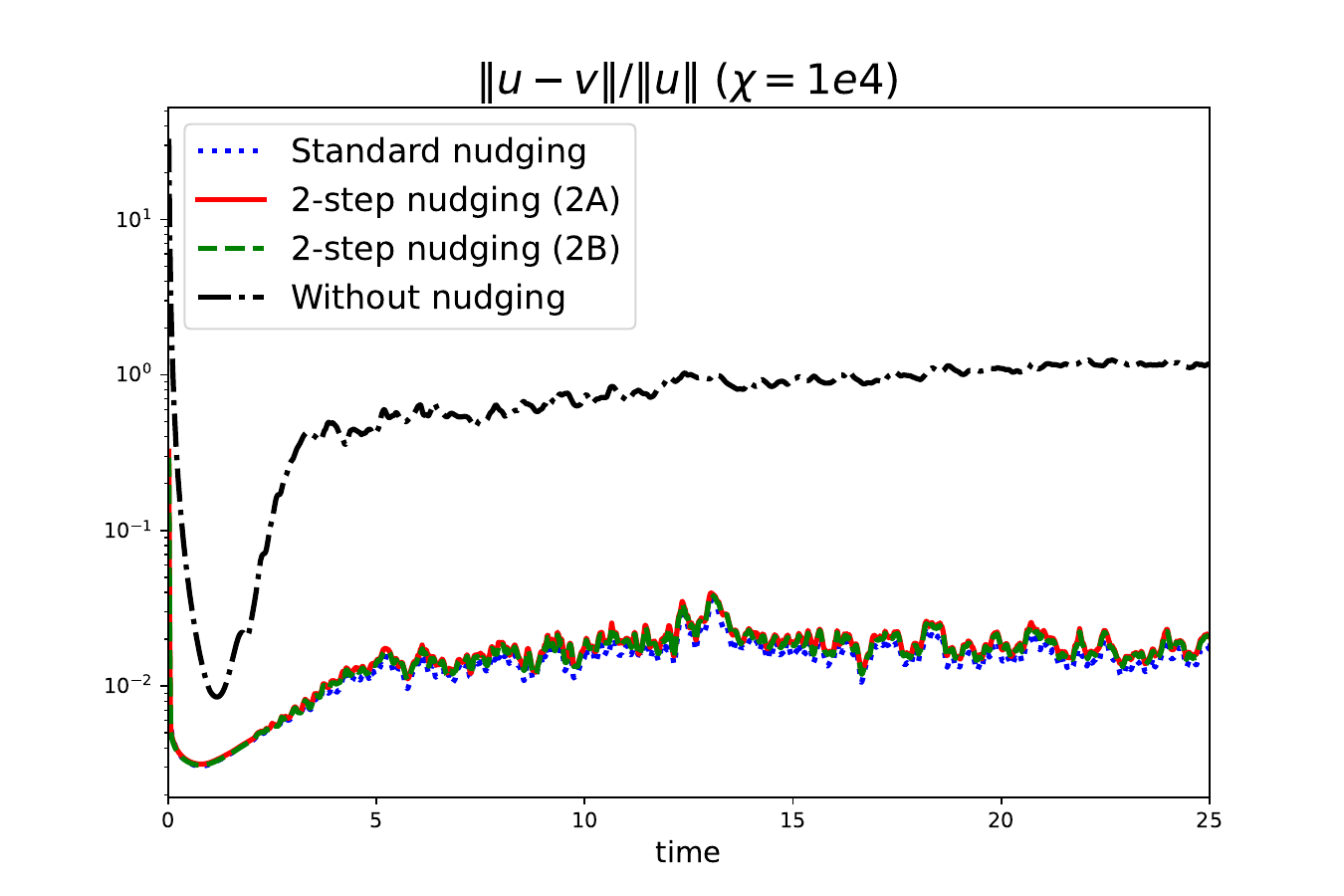}
    \caption{$\chi=1e4$. 2-step nudging methods and standard nudging significantly reduce the velocity error compared to the system without data assimilation.}
    \label{fig:chi_1e4}
\end{figure}

We choose $I_H$ to be the $L^2$ projection and set $\chi = 10^4$. In Figure \ref{fig:chi_1e4}, we compare the relative velocity error of the two-step nudging methods (\ref{eq:Nudging2Step}, \ref{2step_B}), and the standard nudging method (\ref{eq:StandardNudged}), and the case without data assimilation. 2-step nudging with 2A and 2B, and the standard nudging achieve similar result in this test. The relative velocity error decreases rapidly at the beginning for both methods. As the flow evolves into a more complex pattern, the error gradually increases and stabilizes after $t = 15$, reaching a steady-state level on the order of $10^{-2}$. Both nudging methods significantly reduce the velocity error compared to the system without data assimilation, whose relative error remains $\mathcal{O}(1)$ over time.

\begin{figure}[H]
    \centering
    \includegraphics[width=0.9\linewidth]{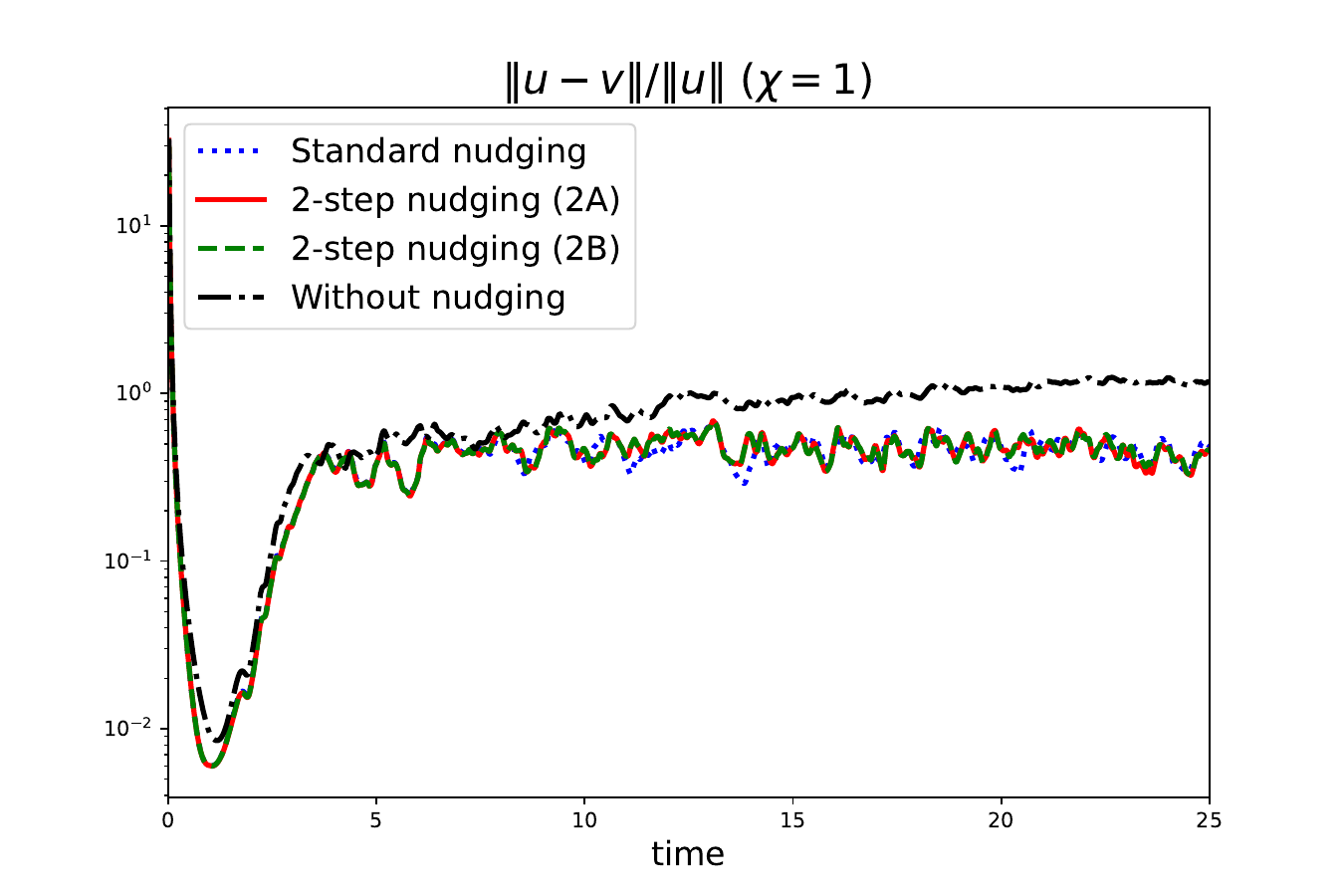}
    \caption{$\chi=1$. We observe a relative velocity error of approximately 0.3 for the 2-step nudging methods with 2A and 2B, and standard nudging method. They improve the $\mathcal{O}(1)$ error without data assimilation.}
    \label{fig:chi_1}
\end{figure}

We also test a much smaller nudging parameter, $\chi = 1$, and observe a relative velocity error of approximately 0.3, as shown in Figure \ref{fig:chi_1}. While this result is worse than with $\chi = 10^4$, it still outperforms the case without data assimilation.

\section{Conclusions}

Section 1.1 summarized several algorithmic difficulties with standard nudging for practical data assimilation. To conclude, we review them and summarize how modular nudging (inspired by the algorithmic form of Kalman filter methods) addresses the difficulties.

\textbf{Legacy codes:} Standard nudging requires an intrusive modification of existing CFD codes. \textit{\ In contrast, modular nudging is minimally intrusive.}

\textbf{Coefficient matrix fill-in:} In standard nudging, if (\ref{eq:StandardNudged}) is discretized on a spatial mesh of size $h\ll H$, the term $\chi I_{H}v^{n+1}$\ couples variables across the observational
macro-elements, increasing coefficient matrix fill-in and operations dramatically. \textit{In contrast, in modular nudging the terms with }$I_{H}$ \textit{ do not affect the coefficient matrix if the explicit update is available and used} 
\begin{equation}\label{explicit_formula}
v^{n+1}=\widetilde {v}^{n+1}-\frac{k\chi }{1+k\chi }I_{H}\widetilde {v}%
^{n+1}+\frac{k\chi }{1+k\chi }I_{H}u(t^{n+1}).
\end{equation}

\textbf{Mesh communication:} Standard nudging also requires communication between the $H-$mesh and the $h-$mesh in assembly of the linear system. Its availability depends upon the software platform used. \textit{In contrast, modular nudging requires only explicit application of prolongation and restriction operators between the
two meshes.}

\textbf{Conditioning:} When $\chi \gg 1$, the nudging term increases the condition number of the associated coefficient matrix. \textit{In contrast,
in modular nudging the parameter }$\chi $\textit{\ does not affect the condition number of any coefficient matrix when (\ref{explicit_formula}) is used.}

\textbf{The }$H,\chi $\textbf{\ conditions:} The precise form of $``H$ 
\textit{small enough} and $\chi $\ \textit{large enough"} in, e.g. 3d are${{{%
{\text { }}}}}\chi \gtrsim {{{{\mathcal{O}}}}}({{{{\mathcal{R}}}}}e^{5}){{{{%
\text { and \ }}}}}H\lesssim {{{{\mathcal{O}}}}}({{{{\mathcal{R}}}}}e^{-3})$%
, smaller than the Kolmogorov micro-scale, \cite{CFLS25}. \textit{In contrast, modular nudging, when these conditions are violated, increases predictability horizons. This new but small result is an open problem for standard nudging and for global in time predictability horizons.}

\textbf{The user-chosen parameter }$\chi $\textbf{\ must be chosen} without analytic guidance. \textit{Recent work \cite{CFLS25} has shown that }$\chi $%
\textit{\ can be chosen self-adaptively with a procedure that can be used for either standard or modular nudging.}

\textbf{Open problems.} The algorithms and analysis herein are only a first step. Many fundamental problems remain, including:

\textit{Analysis of Step 2A with I}$_{H}$\textit{\ an L}$^{2}$\textit{ projection} is an open problem. The issue here is to find an a priori estimate on $\|\nabla v\|$ (and $\|\nabla e\|$) to control error growth caused by vortex stretching.

\textit{Analysis of the impact of standard nudging upon predictability horizons} when the $H$ small and $\chi $\ large conditions fail is an open problem. It must surely extend predictability horizons but it is not clear
how to prove this.

\textit{Numerical analysis of modular nudging for better time discretizations} is an important open problem. Extension to more complex and useful ones, especially time-accurate variable-step, variable-order methods is essential.

To summarize (in the spirit of Richardson):
\begin{quote}\textit{In CFD where nudging can sway,\\ Modularity yields a new, clearer way.\\
When Step 2 ensues,\\
Complexity is reduced,\\
And step-by-step errors decay.\\}
\end{quote}

\end{document}